\documentclass[11pt]{amsart}
\usepackage{amsmath,amssymb,amsthm,bm,enumerate,multido,pstricks,pst-node,caption}
\usepackage[mathscr]{euscript}
\usepackage[foot]{amsaddr}

\newtheorem{theorem}{Theorem}[section]

\newtheorem{lemma}[theorem]{Lemma}
\newtheorem{proposition}[theorem]{Proposition}
\newtheorem{corollary}[theorem]{Corollary}

\theoremstyle{definition}
\newtheorem{definition}[theorem]{Definition}
\newtheorem{example}[theorem]{Example}
\newtheorem{remark}[theorem]{Remark}
\newtheorem{question}{Question}
\newtheorem{construction}[theorem]{Construction}
\def\D{{\mathscr{D}}}
\def\PP{{\mathscr{P}}}

\def\B{{\mathscr{B}}}
\def\C{{\mathscr{C}}}
\def\Aut{{\mathrm{Aut}}}
\def\Sym{{\mathrm{Sym}}}

\def\PG{{\rm PG}}

\def\Ch{\bm{\mathscr{C}}}

\definecolor{amber}{rgb}{1.0, 0.75, 0.0}

\usepackage{anysize}
\marginsize{2cm}{2cm}{2cm}{2cm}

\begin{document}

\title{Block-transitive $2$-designs with a chain of imprimitive point-partitions}
\author{Carmen Amarra$^{\rm\lowercase{a,b,\ast}}$} \author{Alice Devillers$^{\rm \lowercase{a}}$} \author{Cheryl E. Praeger$^{\rm \lowercase{b}}$}
\address{$^{\rm\lowercase{a}}$Centre for the Mathematics of Symmetry and Computation, The University of Western Australia, 35 Stirling Highway, Crawley, Western Australia 6009, Australia}
\email{carmen.amarra@uwa.edu.au, alice.devillers@uwa.edu.au, cheryl.praeger@uwa.edu.au}
\address{$^{\rm\lowercase{b}}${\rm{\emph{P\lowercase{ermanent address}}}}: Institute of Mathematics, University of the Philippines Diliman, C.P. Garcia Avenue, Quezon City 1101, Philippines}
\email{mcamarra@math.upd.edu.ph}
\address{$^{\ast}$Corresponding author}

\maketitle

\begin{abstract}
More than $30$ years ago, Delandtsheer and Doyen showed that the automorphism group of a block-transitive $2$-design, with blocks of size $k$, could leave invariant a nontrivial point-partition, but only if the number of points was bounded in terms of $k$. Since then examples have been found where there are two nontrivial point partitions, either forming a chain of partitions, or forming a grid structure on the point set. We show, by construction of infinite families of designs, that there is no limit on the length of a chain of invariant point partitions for a block-transitive $2$-design. We introduce the notion of an `array' of a set of points which describes how the set interacts with parts of the various partitions, and we obtain necessary and sufficient conditions in terms of the `array' of a point set, relative to a partition chain, for it to be a block of such a design.   

\noindent{This work forms part of the Australian Research Council Discovery Grant project  
DP200100080. }
\end{abstract}

\section{Introduction}

A $t$-$(v,k,\lambda)$ \emph{design} $\D$ is a pair $(\PP,\B)$ with $\PP$ a set of \emph{points} and $\B$ a set of \emph{blocks}, such that $\PP$ has size $v$, each block is a $k$-element subset of $\PP$, and every set of $t$ distinct points is contained in exactly $\lambda$ blocks. A \emph{flag} of $\D$ is a pair $(p,B) \in \PP \times \B$ such that $p \in B$. An \emph{automorphism} of $\D$ is a permutation $g$ of the point set $\PP$ which preserves the block set $\B$. The   design $\D$ is \emph{point-transitive} (respectively, \emph{block-transitive} or \emph{flag-transitive}) if it admits a group $G$ of automorphisms which is transitive on $\PP$ (respectively, $\B$ or the set of all flags). It was proved in 1967 by Block~\cite{Block67} that, if $t\geq2$, then a block-transitive group $G$ is also transitive on the point set. Earlier, in 1961, it had been proved by Higman and McLaughlin \cite{HigMcL61} that, if $t\geq 2$ and $\lambda=1$, then a flag-transitive group $G$ satisfies an even stronger property in its point-action, namely it is point-primitive, which we now define. 

For a point-transitive subgroup $G$ of $\Aut(\D)$, a partition $\C$ of $\PP$ is $G$-invariant if $C^g \in \C$ for any class $C \in \C$ and $g \in G$. The design $\D$ is \emph{$G$-point-primitive} if the only $G$-invariant partitions are the trivial partitions, namely the partition consisting of singletons (which we denote by $\binom{\PP}{1}$) and $\{\PP\}$. The group $G$ is \emph{imprimitive} otherwise. Hence if $G$ is imprimitive then there is a $G$-invariant partition of $\PP$ with $d$ classes each of size $c$, for some $c,d > 1$.

A flag-transitive group is clearly both block-transitive and point-transitive, and the result of  Higman and McLaughlin  \cite{HigMcL61} raised the question of whether flag-transitive groups on $2$-designs were `usually', or `typically' also point-primitive. In 1987, Davies  \cite{D87} showed that this is nearly always the case by proving that, for a fixed $\lambda$, there are only a finite number of possible parameters $v,k$ for which a flag-transitive group on a $2$-$(v,k,\lambda)$ design can be imprimitive on points. At around the same time, in 1989, Delandtsheer and Doyen \cite{DD} proved that nearly all block-transitive $2$-designs are point-primitive, in the sense that, for a fixed $k$, they determined an explicit quantity $f(k)$ such that a block-transitive group on a $2$-design with more than $f(k)$ points must also be point-primitive. On the other hand, Cameron and the third author \cite{CP93} gave several general constructions for families of block-transitive $t$-designs (with $t=2, 3$) which are point-imprimitive. These constructions usually feature a single nontrivial invariant partition, although in some constructions two nontrivial partitions arise, see \cite{grids22, gridsconstr21, CP93, PT03, ZZ2}.

The purpose of this paper is to demonstrate that there is no limit to the `depth of imprimitivity' a block-transitive $2$-design may possess.  We define a family of point-block incidence structures admitting a certain block-transitive group which leaves invariant each partition in a given chain of partitions of the point-set. We call such point actions \emph{$s$-chain-imprimitive}. We  obtain necessary and sufficient conditions for this incidence structure to be a $2$-design (Theorem~\ref{mainthm:2des}). In Section \ref{sec:ex} we give an explicit construction of an infinite family of block-transitive $s$-chain-imprimitive $2$-designs (Construction~\ref{con:s>2}) where the parameter $s$ can be arbitrarily large. Although it took us some time to find the infinite family in Construction~\ref{con:s>2}, we do not know how unusual this construction is.

\begin{question}
    Are there many other infinite families of block-transitive, $s$-chain-imprimitive $2$-designs, with unboundedly large values of $s$?
\end{question}

We show in Theorem~\ref{thm:s>1} that none of the designs in Construction \ref{con:s>2} is a $3$-design and none is flag-transitive. We shall explore the flag-transitive case in an upcoming paper \cite{flagtr-chain}. Moreover we do not know any block-transitive, $s$-chain-imprimitive $3$-designs with $s \geq 3$.

\begin{question}
    For which $s \geq 3$ does there exist an $s$-chain-imprimitive block-transitive $3$-design? In particular, is there a bound on $s$ for such $3$-designs?
\end{question}

\subsection{Statements of the main results}\label{sub:results}

A partition $\C'$ of a set $\PP$ is a \emph{proper refinement} of a partition $\C$ of $\PP$ if each part in $\C'$ is a subset of some part in $\C$, and at least one of these inclusions is proper. We write $\C' \prec \C$ if $\C'$ is a proper refinement of $\C$. In the designs that we consider, a point-transitive subgroup of automorphisms leaves invariant a chain of partitions ordered by $\prec$. We introduce the following concepts.

\begin{definition} \label{def:chain} 
Let $s\geq 2$, let $\PP$ be a set,  and let $G\leq \Sym(\PP)$. 
We say that the action of $G$ on $\PP$ is \emph{$s$-chain-imprimitive} if $G$ is transitive and
there exist $G$-invariant partitions $\C_0, \C_1, \ldots, \C_s$ of $\PP$, which form a chain
    \begin{equation} \label{eq:chain}
    \Ch : \ {\textstyle\binom{\PP}{1}} = \C_0 \prec \C_1 \prec \cdots \prec \C_{s-1} \prec \C_s = \{\PP\}. 
    \end{equation}
under the partial order $\prec$. If $\D$ is a design with point set $\PP$, such that $G \leq \Aut(\D)$ acts $s$-chain-imprimitively on  $\PP$, then we say that $\D$ is \emph{$(G,s)$-chain-imprimitive}.
\end{definition}

Thus a $2$-chain-imprimitive design is a point-imprimitive design with a specified nontrivial invariant point-partition. In Subsection~\ref{sub:history} we mention some results in the literature on block-transitive, point-imprimitive $2$-designs to provide further context for our work. Less is known about $(G,s)$-chain-imprimitive designs with $s\geq3$, although we can make geometric constructions exploiting certain number theoretic properties. For example, if $\D$ is a Desarguesian projective plane of order $q$ such that the number $q^2+q+1$ of points has at least $s$ pairwise distinct prime factors, and if $G = C_{q^2+q+1}$ is a `Singer cycle' acting regularly on the points, and on the lines, of $\D$ then $\D$ is $(G,s)$-chain-imprimitive and $G$-block-transitive. However we found, after some computer experimentation, that it was not easy to find large families of prime powers $q$ for which $q^2+q+1$ has many distinct prime factors. Thus the first, rather natural, question we set for ourselves was the following.

\begin{center} 
\emph{For which values of $s\geq 3$, can we give explicit constructions for\\ infinite families of block-transitive, $s$-chain-imprimitive $2$-designs?} 
\end{center}

With a lot of effort we were able to give a complete answer to this question, namely the answer is ``all values of $s$'', via Constructions \ref{con:s=2} (for $s=2$), \ref{con:s=3} (for $s=3$) and \ref{con:s>2} (for arbitrary $s$). Thus we have:

\begin{theorem} \label{thm:ex}
Given any integer $s \geq 2$, there exist infinitely many $2$-designs which are $G$-block-transitive and $(G,s)$-chain-imprimitive for some automorphism group $G$.
\end{theorem}

The designs produced by Constructions \ref{con:s=3} and \ref{con:s>2} are generalisations of the Cameron-Praeger designs described in \cite[Section 2]{CP93}. A Cameron-Praeger design $\D = (\PP,\B)$ has automorphism group $G$ which is the full stabiliser in $\Sym(\PP)$ of a specified nontrivial point-partition $\C$, and block set $\B=B^G=\{B^g\mid g\in G\}$ for some $B\subseteq \PP$. In \cite{CP93} Cameron and the third author obtained necessary and sufficient conditions on certain parameters associated with $B$ and $\C$ for $\D$ to be a $2$-design, and to be a $3$-design, see \cite[Proposition 2.2]{CP93}. These parameters are the integers $x_C := |B \cap C|$ for $C \in \C$, which describe how $B$ is partitioned among the classes of $\C$.  In Section~\ref{subsec:array} we develop this `$1$-dimensional array', formed by the $x_C$ to describe $B$ relative to a single nontrivial partition $\C$, into an `array function' which describes a point subset $B$ relative to an $s$-chain of point partitions, and encodes the size $x_C=|B\cap C|$ for any part $C$ of any of the nontrivial partitions in the $s$-chain (Definition~\ref{def:array} and \eqref{array}). We introduce a notion of equivalence on these array functions (Definition~\ref{def:equivarray}) and show that equivalent array functions correspond precisely to orbits $B^G$ of point subsets under the iterated wreath product $G$ (Lemma~\ref{lem:array-equiv}). Properties of these array functions are critical for studying $s$-chain-imprimitive $2$-designs, and in particular for proving   Theorem \ref{mainthm:2des}, which generalises \cite[Proposition 2.2]{CP93} for $s$-chain-imprimitive $2$-designs.  Theorem \ref{mainthm:2des} is used to verify in the proof of Theorem~\ref{thm:s>1} that Construction~\ref{con:s>2} does indeed produce $2$-designs.

\begin{theorem} \label{mainthm:2des}
Given integers $s \geq 2$ and $e_1, \ldots, e_s \geq 2$. Let $G = S_{e_1} \wr \ldots \wr S_{e_s}$, acting $s$-chain-imprimitively on a set $\PP$ of size $v = \prod_{i=1}^s e_i$ and with invariant partitions which form an $s$-chain $\Ch$ under $\prec$ as in \eqref{eq:chain}, such that for $1 \leq i \leq s$ each class of $\C_i$ contains $e_i$ classes of $\C_{i-1}$. Let $\D = (\PP,\B)$ where $\B = B^G$ for some $k$-subset $B$ of $\PP$. Then $\D$ is a $G$-block-transitive, $(G,s)$-chain-imprimitive $2$-design if and only if the following equations hold:
    \begin{equation} \label{2des-1}
    \sum_{C \in \C_1} x_{C} \left( x_{C} - 1 \right)
    = \frac{k(k-1)}{v-1} (e_1 - 1)
    \end{equation}
and
    \begin{equation} \label{2des}
    \sum_{C \in \C_{i-1}} x_C \left( x_{C^+} - x_C \right) = \frac{k(k-1)}{v-1} (e_i - 1) \prod_{j \leq i-1} e_j \quad \text{for $i \in \{2, \ldots, s\}$},
    \end{equation}
where $C^+$ denotes the unique $\C_i$-class which contains $C$, and $x_C := |B \cap C|$ and $x_{C^+} := | B \cap C^+|$. Furthermore if \eqref{2des-1} and \eqref{2des} are satisfied, then the quantities in each of the equalities \eqref{2des-1} and \eqref{2des} are even integers, and $v-1$ divides $\binom{k}{2} \cdot \gcd(e_1-1,\dots, e_s-1)$.
\end{theorem}

\begin{remark} \label{rem:s-cond}
Equalities \eqref{2des-1} and \eqref{2des} consist of $s$ conditions, and each of these corresponds to an $i \in \{1, \ldots, s\}$. In particular, condition \eqref{2des-1} is \eqref{2des} for $i=1$ expressed in a ``nicer" way.
\end{remark}

The group $S_{e_1} \wr \ldots \wr S_{e_s}$ in Theorem \ref{mainthm:2des} is the full stabiliser of the chain $\Ch$ of partitions. The following result, for the more general case where the block-transitive automorphism group $G$ is an arbitrary $s$-chain-imprimitive subgroup of $S_{e_1} \wr \ldots \wr S_{e_s}$,  is a consequence of Theorem \ref{mainthm:2des} and \cite[Proposition 1.1]{CP93}.

\begin{theorem} \label{thm:subdes}
Given integers $s \geq 2$ and $e_1, \ldots, e_s \geq 2$. Let $\PP$ be a set of size $v = \prod_{i=1}^s e_i$, with an $s$-chain $\Ch$, as in \eqref{eq:chain}, of partitions of $\PP$, such that for each $1 \leq i \leq s$ each class of $\C_i$ contains $e_i$ classes of $\C_{i-1}$. Let $\D = (\PP,\B)$, where $\B = B^G$ for some $k$-subset $B$ of $\PP$ and some $G \leq S_{e_1} \wr \ldots \wr S_{e_s}$ that leaves invariant the chain $\Ch$. If $\D$ is a $G$-block-transitive, $(G,s)$-chain-imprimitive $2$-design then conditions \eqref{2des-1} and \eqref{2des} hold.
\end{theorem}

\begin{remark}
For the designs $\D$ in Theorem \ref{thm:subdes}, the number $b$ of blocks is given by $b = |G:G_B|$, where $G_B$ is the setwise stabiliser of $B$, and from this we can obtain $\lambda = \frac{bk(k-1)}{v(v-1)}$. We perform these computations explicitly for the designs in Construction \ref{con:s>2}, as well as for the examples in Section \ref{sec:uniform}. We note that for all of these the values of $\lambda$ are rather large, with $\lambda$ greater than the size of a $\C_{s-1}$-class.
\end{remark}

In the original paper of Delandtsheer and Doyen~\cite{DD} on block-transitive, point-imprimitive $2$-designs $\D$, the notions of inner and outer point-pairs were introduced to describe pairs of points in a block of $\D$ which did, or did not, lie together in the same class of a specified nontrivial point partition. In Section~\ref{s:innout}, we extend these notions to describe point-pairs in the context of $s$-chains of partitions, which proves helpful in our exposition. 

In Section \ref{sec:uniform} we give examples of block-transitive, $s$-chain-imprimitive $2$-designs that do not fall under the constructions in Section \ref{sec:ex} mentioned above. For these examples the underlying $s$-chain of point-partitions is \emph{uniform}, that is,  the cardinality $|\C_{i-1}(C)|$ is constant, independent of $i$ and $C\in\C_i$. We found several examples with $s=3$ and one example with $s=4$ (see Examples~\ref{ex:e=3} and~\ref{ex:s=4}), but we were unable to generalise these to give an infinite family. We note that none of the $s$-chains of partitions in the examples in Construction \ref{con:s>2} are uniform.

\begin{question}\label{Q4}
Are there infinitely many block-transitive, $s$-chain-imprimitive $2$-designs relative to a uniform $s$-chain of partitions with $s\geq3$? In particular, are there any examples with $s\geq5$?
\end{question}

\subsection{Commentary on block-transitive, point-imprimitive designs}\label{sub:history}

Delandtsheer and Doyen showed in \cite{DD} that, if a $2-(v,k,\lambda)$ design is block-transitive, then $v \leq \left(\binom{k}{2} - 1\right)^2$. In \cite{CP93}, Cameron and the third author studied point-imprimitive $t$-designs for arbitrary $t$, and showed that if such a design  is block-transitive then $t\leq 3$, and if it is flag-transitive then $t\leq 2$. They also showed, for a flag-transitive, point-imprimitive $2-(v,k,\lambda)$ design, that $v \leq (k-2)^2$, and they conjectured that $3-(v,k,\lambda)$ designs should also be further restricted: namely if such a $3$-design is block-transitive and point-imprimitive, then $v \leq \binom{k}{2} + 1$ should hold. Mann and Tuan in \cite{MT03} proved this conjecture, and Tuan in \cite{Tuan} constructed infinitely many examples of $3$-designs attaining this bound. Recent work has seen classification of $3$-designs, and of flag-transitive $2$-designs, with small parameters if the point action is imprimitive: in \cite[Theorem 2]{ZPW} for block-transitive $3-(v,k,\lambda)$ designs with  $k\leq 6$ (all have  $(v,k)=(16,6)$); of flag-transitive symmetric $2$-designs with $\lambda\leq 4$ in \cite{LPR, OR, P, PZ}, leading up to a classification for $\lambda\leq 10$ in \cite{MS, M, M2}; and of flag-transitive non-symmetric $2$-designs with $\lambda\leq 4$ in  \cite{DLPX, DP21, DP22}. Imprimitive flag-transitive designs with other small values of $\lambda$ have been studied in \cite{ZZ}.

There have been studies focusing on point-imprimitive designs with two distinct nontrivial invariant partitions which do not form a chain. 
Here the point set has size $v = cd$ (with $c,d > 1$) and the automorphism group $G$ preserves two partitions $\C_1$ and $\C_2$, one with $d$ classes of size $c$, the other with $c$ classes of size $d$, and such that $|C_1 \cap C_2| = 1$ for any $C_1 \in \C_1$ and $C_2 \in \C_2$. Let us call such designs \emph{grid-imprimitive} since the points of $\PP$ can be identified with the `grid'  $\{1, \ldots, c\} \times \{1, \ldots, d\}$ so that the invariant partitions are the set of rows and the set of columns. In 1993, Cameron and the third author gave a general construction for block-transitive, grid-imprimitive  $t$-designs for $t=2, 3$, where the point set was reinterpreted as the edge-set of a complete bipartite graph, so that blocks were certain edge-induced subgraphs, see \cite[Proposition 3.6]{CP93}. Although the necessary and sufficient conditions for a $2$-design in \cite[Proposition 3.6]{CP93} are purely  arithmetic, one of the conditions for a $3$-design exploited the graph theoretic interpretation, namely a restriction on the number of paths of length $3$ in the subgraphs corresponding to blocks. A fully graph theoretic approach was recently developed by Alavi, Daneshkhah, and the second and third authors in \cite{grids22} giving necessary and sufficient conditions for grid-imprimitive  $2$-designs and $3$-designs in graph theoretic terms. The characterisation led to new explicit constructions of these designs in \cite{grids22}, and also allowed geometrical recognition of some of the designs in the classifications by Zhan et.al. in \cite[Theorems 1 and 2]{ZZBPG} of block-transitive  $2-(v,4,\lambda)$ designs (achieved largely by computer search using \textsc{Magma} \cite{magma}). 
An alternative, more combinatorial approach was used in \cite{gridsconstr21} by Brai\'{c} et.al. to develop the ideas in \cite{CP93}, and this was exploited, again using \textsc{Magma} \cite{magma}, to construct all flag-transitive grid-imprimitive $2$-designs with automorphism group $S_c \times S_d$, where $4 \leq c \leq d \leq 70$.
In the context of linear spaces, that is, $2$-$(v,k,1)$ designs, it was remarked in \cite[p.~1402]{Linsp09} that ``for each parameter set for which an imprimitive line-transitive linear space is known to exist, there is an example where the line-transitive, point-imprimitive group G preserves a grid structure on the point set''. For example, (up to isomorphism) exactly $27$ of the $467$ point-imprimtive block-transitive $2-(729,8,1)$ designs are grid-imprimitive, see \cite{NNOPP}. Again for these elusive structures it is possible to deal with the case where the parameters are small:  the classification begun in \cite{ Linsp09grid} of block-transitive grid-imprimitive $2$-$(v,k,1)$ designs with $\gcd(k,v-1) \leq 8$ was completed in \cite{Linsp09}, where the problem \cite[Problem 2]{Linsp09} of finding more block-transitive grid-imprimitive linear spaces was posed. 

\emph{Acknowledgement.} The authors wish to thank their colleagues in the Centre for the Mathematics of Symmetry and Computation of The University of Western Australia for the opportunity to begin work on this project at the annual CMSC research retreat in February 2019. The authors are grateful to the anonymous referees for their comments, in particular the suggestion to determine the values of $\lambda$.

\section{Chain structures and iterated wreath products} \label{sec:2-des}

\subsection{Chain structure on a set} \label{subsec:chain}

Let $s \geq 2$ and let $e_1, \ldots, e_s$ be integers with $e_i \geq 2$ for each $i$. Let $\PP$ be a set of size $|\PP| = \prod_{i=1}^s e_i$, and for each $i \in \{1, \ldots, s\}$ let $\Delta_i$ be a set of size $|\Delta_i| = e_i$. Then we can identify  $\PP$ with the set of $s$-tuples in $\prod_{i=1}^s \Delta_i$.

For any $i \in \{0, \ldots, s-1\}$ and any fixed $(s-i)$-tuple $(\delta_j)_{j>i} = (\delta_{i+1},\delta_{i+2},\ldots,\delta_{s-1},\delta_{s})
\in \prod_{j>i} \Delta_j$, let $C_{(\delta_j)_{j>i}}$ denote the set
    \begin{equation} \label{class}
    C_{(\delta_j)_{j>i}} := \big\{ (\varepsilon_j)_{j=1}^s \in \PP \ \big| \ \varepsilon_j = \delta_j \ \text{for } \ i+1 \leq j \leq s \big\}.
    \end{equation}
Note in particular that $C_{(\delta_j)_{j>0}}$ is the singleton set $\big\{ (\delta_j)_{j>0} \big\}$. Also by convention we write
    \[ C_{(\delta_j)_{j>s}} := \PP, \]
and we consider $\prod_{j>s} \Delta_j$ as a singleton. For $i \in \{0, \ldots, s\}$ define the family $\C_i$ of subsets of $\PP$ as
    \begin{equation} \label{partition}
    \C_i := \Big\{ C_{(\delta_j)_{j>i}} \ \Big| \ (\delta_j)_{j>i} \in \prod_{j>i} \Delta_j \Big\}.
    \end{equation}
Then $\C_0 = \binom{\PP}{1}$, the set of singleton subsets of $\PP$; $\C_s = \{\PP\}$; and for $1 \leq i \leq s-1$ the family $\C_i$ is a nontrivial partition of $\PP$. For $0 \leq i \leq s-1$ the partition $\C_i$ is a proper refinement of $\C_{i+1}$ since each $e_j \geq 2$, and we write $\C_i \prec \C_{i+1}$ to denote this 
 relation. Each class of $\C_{i+1}$ contains $e_{i+1}$ classes of $\C_i$. Thus we have
    \[ \C_0 \prec \C_1 \prec \cdots \prec \C_{s-1} \prec \C_s, \]
so the partitions $\C_0, \ldots, \C_s$ form an $s$-chain $\Ch$ under the relation (partial order) $\prec$.

Set $e_0 := 1$ and $\Delta_0$ to be a singleton, and by convention take $\prod_{j > i} e_j := 1$ for $i = s$. Then for each $i \in \{0, \ldots, s\}$ the partition $\C_i$ consists of $d_i$ classes each of size $c_i$, where
    \begin{equation} \label{cd}
    c_i := \prod_{j \leq i} |\Delta_j| = \prod_{j \leq i} e_j \quad \text{and} \quad d_i := \prod_{j > i} |\Delta_j| = \prod_{j > i} e_j.
    \end{equation}
In particular $c_0 = d_s = 1$, $c_s = d_0 = |\PP|$, $c_1 = e_1$, and $d_{s-1} = e_s$.

\def\chainclass{
    \multido{\n=-0.5+1.0}{2}{\psline(0,\n)(1.25,\n)}
    \multido{\n=1.25+1.25}{3}{\psline[linewidth=1pt](\n,-1.5)(\n,1.5)}
    \pspolygon[linewidth=1.5pt](5,1.5)(0,1.5)(0,-1.5)(5,-1.5)
    }
\begin{figure}
    \centering
    \begin{pspicture}(0,-3)(5,3)
    \chainclass
    \pnode(1.25,1.5){a} \pnode(2.5,1.5){b} \ncbar[nodesep=4pt,angle=90,armA=5pt,armB=5pt]{a}{b} \rput(1.875,2.25){$e_{i+1}$ columns ($\C_i$-classes)}
    \pnode(0,-0.5){c} \pnode(0,0.5){d} \ncbar[nodesep=4pt,armA=5pt,armB=5pt,angle=180]{c}{d} \rput(-2,0){\parbox[c]{2.5cm}{\flushright $e_i$ small \\ rectangles \\ ($\C_{i-1}$-classes) \\ per column}}
    \pnode(0,-1.5){e} \pnode(5,-1.5){f} \ncbar[nodesep=4pt,armA=5pt,armB=5pt,angle=-90]{e}{f} \rput(2.5,-2.25){$\C_{i+1}$-class (big rectangle)}
    \end{pspicture}
    \caption{Visual interpretation of a $\C_{i+1}$-class, where $1 \leq i < s$}
    \label{fig:class}
\end{figure}

\subsection{Inner and outer point-pairs}\label{s:innout}

Delandtsheer and Doyen \cite{DD} introduced the terms \emph{inner pair} and \emph{outer pair}  to refer to unordered pairs of distinct points with a specific relationship with respect to a nontrivial partition $\C$ of $\PP$ (two points form an inner pair if they lie in the same class of $\C$, and an outer pair if they lie in different classes of $\C$). We generalise this terminology to pairs of points with respect to an $s$-chain $\Ch$ of nontrivial partitions, as in the above. A pair $\{\alpha,\beta\}$ of distinct points is called an \emph{$i$-outer pair} (for $0 \leq i \leq s-1$) if $\alpha$ and $\beta$ lie in different $\C_i$-classes; an \emph{$i$-inner pair} (for $1 \leq i \leq s$) if $\alpha$ and $\beta$ lie in the same $\C_i$-class; and an \emph{$i$-strong inner pair} (for $1 \leq i \leq s$) if it is an $i$-inner pair and an $(i-1)$-outer pair. Of particular significance in later proofs is the set of all $i$-strong inner pairs for each $i \in \{1, \ldots, s\}$, which we denote by $O_i$, that is,
    \begin{equation} \label{O_i}
    O_i = \left\{ \{\alpha,\beta\} \ \vline \ \text{$\alpha$ and $\beta$ belong to the same $\C_i$-class but to different $\C_{i-1}$-classes} \right\}.
    \end{equation}
For each $i \in \{1, \ldots, s\}$, the set $O_i$ has size
    \begin{equation} \label{|O_i|}
    |O_i| = \frac{1}{2}\,|\PP| \cdot \left( \prod_{j \leq i} e_j - \prod_{j \leq i-1} e_j \right) = \frac{1}{2}\, |\PP| \cdot (e_i - 1) \prod_{j \leq i-1} e_j,
    \end{equation}
and in particular
    \[ |O_1| = \frac{1}{2}\, |\PP| \cdot (e_1 - 1). \]

\subsection{Iterated wreath product groups}

Let $s$, $\Delta_i$, $e_i$, and $\C_i$ be as in Subsection \ref{subsec:chain}, and for each $i \in \{1, \ldots, s\}$ let $G_i$ be a subgroup of the symmetric group $ \Sym(\Delta_i)\cong S_{e_i} $. Also, following \cite{BPRS}, for each $i \in \{1, \ldots, s\}$ let 
    \begin{equation} \label{F_i}
    F_i := \Big\{ f \ \big| \ f \ \text{is a function from $\textstyle\prod_{j>i} \Delta_j$ to $G_i$} \Big\},
    \end{equation}
which is a group under pointwise multiplication defined by
    \[ ((\delta_j)_{j>i}) (f \cdot h) := ((\delta_j)_{j>i})f \cdot ((\delta_j)_{j>i})h \]
for all $f, h \in F_i$ and $(\delta_j)_{j>i} \in \prod_{j>i} \Delta_j$. Let
    \begin{equation} \label{wreath}
    G := \prod_{i=1}^s F_i,
    \end{equation}
and for each $(f_1, \ldots, f_s) \in G$ define the map $\widehat{f}_{j>i} : \prod_{j>i} \Delta_j \rightarrow \prod_{j>i} \Delta_j$ by
    \begin{equation} \label{induced}
    ((\delta_j)_{j>i}) \widehat{f}_{j>i} = \Big( \delta_j^{((\delta_\ell)_{\ell > j})f_j} \Big)_{j>i},
    \end{equation}
and where for $j<s$, $(\delta_\ell)_{\ell > j}$ denotes the unique point of $\prod_{\ell>j} \Delta_\ell$ determined by $(\delta_j)_{j>i}$, and for $j=s$,  $(\delta_\ell)_{\ell>j}$ denotes the unique element of the empty product $\prod_{\ell>s} \Delta_\ell$. It is shown in \cite{BPRS} (see Lemma 4 and Theorem A) that $G$ is a group under the operation
    \begin{equation} \label{operation}
    fh := \big( f_i \cdot (\widehat{f}_{j>i} h_i) \big)_{i=1}^s \ \text{for all} \  f= (f_1, \ldots, f_s), \,h = (h_1, \ldots, h_s)  \in G,
    \end{equation}
where $\widehat{f}_{j>i} h_i$ denotes the composition of $\widehat{f}_{j>i}$ and $h_i$, and $f_i \cdot (\widehat{f}_{j>i} h_i)$ denotes the pointwise product of $f_i$ and $\widehat{f}_{j>i} h_i$. The group $G$ is the \emph{iterated wreath product} of the groups $G_1, \ldots, G_s$.

The iterated wreath product $G$ acts on $\PP$ via the map which sends any pair $\big((\delta_1, \ldots, \delta_s), \,f\big) \in \PP \times G$ to $(\delta_1, \ldots, \delta_s)^f \in \PP$, where
    \begin{align} \label{action}
    &\left( \delta_1, \ldots, \delta_s \right)^f := \left( \varepsilon_1, \ldots, \varepsilon_s \right), \quad \ \text{with} \ \varepsilon_i = \delta_i^{ \left( (\delta_j)_{j>i} \right)f_i} \ \text{for} \ 1 \leq i \leq s.
    \end{align}
The action \eqref{action} is a faithful group action, so we may identify $G$ (via this action) with  a permutation group on $\PP$. Observe that for $s = 2$ the group $G = G_1 \wr G_2$ is a simple wreath product and the action \eqref{action} is equivalent to the standard imprimitive action of $G_1 \wr G_2$ on $\Delta_1 \times \Delta_2$.

Recall from \eqref{class} and \eqref{partition} that each class of the partition $\C_i$ consists of points that agree in their $j$th coordinates for all $j > i$. Hence, for any fixed $f \in G$, and for any two points $\gamma$ and $\delta$ that belong to the same $\C_i$-class $C$, the images $\gamma^f$ and $\delta^f$ under the map \eqref{action} are again points that belong to the same $\C_i$-class $C'$ (possibly different from $C$). So for each $i \in \{1, \ldots, s-1\}$ the group $G$ leaves invariant the nontrivial partition $\C_i$ (as defined in \eqref{partition}), and induces an action on $\C_i$ given by
    \[ 
    \big( C_{(\delta_j)_{j>i}} \big)^f = C_{((\delta_j)_{j>i}){\widehat{f}_{j>i}}} 
    \]
for every $(\delta_j)_{j>i} \in \prod_{j>i} \Delta_j$ and $f = (f_1, \ldots, f_s) \in G$, and where $\widehat{f}_{j>i}$ is as defined in \eqref{induced}. The kernel of the $G$-action on $\C_i$ is $G_{(\C_i)} \cong (G_1 \wr \ldots \wr G_i)^{\prod_{j>i} e_j}$ and the group induced by $G$ on $\C_i$ is $G^{\C_i} = G/G_{(\C_i)}  \cong G_{i+1} \wr \ldots \wr G_s$.

The 1983 paper \cite{BPRS} by Bailey, Praeger, Rowley, and Speed considered general posets of partitions of a given set, of which the chain of partitions described above is a special case. Their result \cite[Theorem B]{BPRS} describes the stabiliser $G$ of such a general poset structure, and \cite[Theorem C]{BPRS} gives the orbits of $G$ on the set $\PP \times \PP$ of ordered pairs. We state below the specialisation of these results for chains of partitions, which is a direct outcome of \cite[Theorems B and C]{BPRS} and the discussion above.

\begin{theorem} \label{thm:BPRSchain}
Let $s \geq 2$ and for each $i \in \{1, \ldots, s\}$ let $e_i \in \mathbb{Z}$ with $e_i \geq 2$, let $\Delta_i$ be a set with $|\Delta_i| = e_i$, let $\PP = \prod_{i=1}^s \Delta_i$, and let $\C_i$ be as in \eqref{partition}. Also let $\C_0$ denote the partition of $\PP$ into singletons. Then the following all hold.
    \begin{enumerate}[(i)]
    \item The partitions $\C_i$ form an $s$-chain $\C_0 \prec \C_1 \prec \cdots \prec \C_s$ (relative to the partial order $\prec$), and its stabiliser in $\Sym(\PP)$ is the subgroup $G$  defined in \eqref{wreath}--\eqref{operation}, where each $G_i = \Sym(\Delta_i)$, and moreover $G$ is the iterated wreath product $G = S_{e_1} \wr \ldots \wr S_{e_s}$.
    \item For $1 \leq i \leq s-1$, the group $G^{\C_i}$ induced by $G$ on the partition $\C_i$ is permutationally isomorphic to the group induced by $G$ on $\prod_{j>i} \Delta_j$, and thus
        $G^{\C_i} \cong S_{e_{i+1}} \wr \cdots \wr S_{e_s}.$
    \item The orbits of $G$ on the set of unordered pairs of elements in $\PP$ are the sets $O_i$, $1 \leq i \leq s$, of $i$-strong inner pairs, as defined in \eqref{O_i}.
    \end{enumerate}
\end{theorem}

\subsection{Array of a subset of $\PP$} \label{subsec:array}

We use the notation of Subsection \ref{subsec:chain}. In the case where $s = 2$ there is only one nontrivial partition $\C = \C_1$ of $\PP$ in the chain $\Ch$. If $\C$ has $d = d_1$ classes each of size $c = c_1$ then the full stabiliser of $\C$ is the wreath product $G = S_c \wr S_d$. Recall that $e_1 = c_1$ and $e_2 = d_{s-1} = d_1$, so $G = S_{e_1} \wr S_{e_2}$. Let $\D = (\PP,\B)$, where $\B = B^G$ for some nonempty subset $B$ of $\PP$ of size $k < |\PP| = cd$. In \cite{CP93}, Cameron and the third author observed that the subset $B$ has an associated $d$-tuple $(x_1, \ldots, x_d)$, where $x_i$ is the number $|B \cap C|$ of points of $B$ that belong in the $i$th part $C$ of $\C$. Hence $0 \leq x_i \leq c$ for each $i$ and $\sum_{i=1}^d x_i = k$. Define $\mathbf{x}$ to be the $d$-tuple $(x'_1, \ldots, x'_d)$ obtained by rearranging the  coordinates of $(x_1, \ldots, x_d)$ so that they occur in non-increasing order. Any subset $B' \in \B$ is partitioned by the classes of $\C$ into subsets with sizes given by the entries of $\mathbf{x}$ in some order, and in fact $\B$ consists of all $k$-subsets of $\PP$ for which the $d$-tuple of class intersection sizes is equal to some permutation of the entries of $\mathbf{x}$. Hence $\D$ is determined by $c$, $d$, and $\mathbf{x}$, and is denoted $\D(c,d;\mathbf{x})$.

We now consider the case of arbitrary $s \geq 2$, and seek a generalisation of the $d$-tuple $\mathbf{x}$ to describe how subsets of $\PP$ intersect the various classes of the partitions $\C_i$. In this case there is more than one nontrivial partition $\C_i$ of $\PP$ when $s\geq 3$. To any nonempty subset $B$ of $\PP$, we associate a function as follows. 
    \begin{equation} \label{array}
    \chi_B : \bigcup_{i=1}^{s} \C_i \rightarrow \mathbb{Z}_{\geq 0}, \quad
    \chi_B(C) := \left| B \cap C \right| \quad
    \text{for each $i \in \{1, \ldots, s\}$ and $C \in \C_i$.}
    \end{equation}
We call the function $\chi_B$ in \eqref{array} the \emph{array function} of the subset $B$ with respect to the chain $\Ch$ of partitions. For example, if $s=2$, then for a subset $B$ of $\PP$, its associated $d$-tuple in the previous paragraph is $(\chi_B(C_1), \dots, \chi_B(C_d))$, where $\C_1=\{C_1,\dots,C_d\}$. (Note that the array of $B$ is determined by the values of $\chi_B$ on the $\C_i$-classes for $1 \leq i \leq s-1$, since for $C \in \C_s$ we have $C = \PP$, so $\chi_B(C) = |B|$.)

More generally we say that a function $\chi : \bigcup_{i=1}^{s} \C_i \rightarrow \mathbb{Z}_{\geq 0}$ is \emph{an array function with respect to $\Ch$} if $\chi = \chi_B$ for some subset $B$ of $\PP$, with $\chi_B$  as in \eqref{array}. Also we make a synthetic definition of an abstract array function which requires only certain conditions on the function values. It uses the following notational convention: for $i \in \{1, \ldots, s-1\}$ and any nonempty subset $X \subseteq \PP$, let
    \begin{equation} \label{C(X)}
    \C_i(X) := \{ C \in \C_i \ | \ X \cap C \neq \varnothing \}.
    \end{equation}
Thus in particular, if $C \in \C_j$ and $i < j$ then $\C_i(C)$ is the set of all $\C_i$-classes contained in $C$, and if $i > j$ then $\C_i(C)$ is the set consisting of the unique $\C_i$-class which contains $C$. If $j = i-1$ we denote this class by $C^+$, that is,
    \begin{equation} \label{C+}
    C^+ := \text{unique $\C_i$-class which contains $C \in \C_{i-1}$}.
    \end{equation}
In particular, if $C = C_{(\delta_\ell)_{\ell>i-1}}$, then $C_{(\delta_\ell)_{\ell>i-1}}^+ = C_{(\delta_\ell)_{\ell>i}}$.

\begin{definition}\label{def:array}
Let $s\geq2$, suppose that $\Ch$ is an $s$-chain as in \eqref{eq:chain} in Definition~\ref{def:chain}, and let $c_i$ be as in \eqref{cd}. Then a function  $\chi : \bigcup_{i=1}^{s} \C_i \rightarrow \mathbb{Z}_{\geq 0}$ is called an \emph{(abstract) array function} with respect to $\Ch$ if and only if the following two conditions hold:
    \begin{enumerate}[(i)]
    \item \label{array1} for each $i \in \{1, \ldots, s\}$ and each $C \in \C_i$, \ $\chi(C) \leq c_i$, and
    \item \label{array2} for each $i \in \{2, \ldots, s\}$ and each $C \in \C_i$, \ $\sum_{C' \in \C_{i-1}(C)} \chi(C') = \chi(C)$.
    \end{enumerate}
\end{definition}

We show that array functions and abstract array functions are equivalent concepts. 

\begin{lemma} \label{lem:array}
Let $s \geq 2$, and let $\Ch$ be an $s$-chain of point partitions as in \eqref{eq:chain} in Definition~$\ref{def:chain}$. Then a function $\chi : \bigcup_{i=1}^{s} \C_i \rightarrow \mathbb{Z}_{\geq 0}$ is an array function with respect to $\Ch$ if and only if $\chi$ is an abstract array function with respect to $\Ch$.
\end{lemma}

\begin{proof}
Assume first that $\chi = \chi_B$ for some $B \subseteq \PP$. Then for any $i \in \{1, \ldots, s-1\}$ and $C \in \C_i$, we have $\chi(C) = |B \cap C| \leq |C|$ by \eqref{array}, and $|C| = c_i$ by \eqref{cd} and the preceding line of text. Thus $\chi(C) \leq c_i$ and \eqref{array1} holds. Also, for any $i \in \{1, \ldots, s\}$ and any $\C_i$-class $C$ we have
    \[ 
    \sum_{C' \in \C_{i-1}(C)} \chi(C') = \sum_{C' \in \C_{i-1}(C)} |B \cap C'| = |B \cap C| = \chi(C), 
    \]
which is condition \eqref{array2} of Definition~\ref{def:array}. Thus $\chi$ is an abstract array function.

Conversely, assume that $\chi$ is an abstract array function, so conditions \eqref{array1} and \eqref{array2} hold. We define a subset $B$ of $\PP$ as follows: by \eqref{array1}, for each $\C_1$-class $C$ we have $\chi(C) \leq c_1 = |C|$, and hence we may choose a subset $B_C$ of $C$ of size $\chi(C)$ (the empty set if $\chi(C)=0$). Define $B:=\bigcup_{C\in\C_1} B_C$. By definition $\chi(C)=|B\cap C|$ for each $C\in\C_1$. We will show that this property holds also for each $\C_i$ with $2 \leq i \leq s$; from this it will follow that $\chi = \chi_B$, by \eqref{array}. Proceeding by induction on $i$, assume that $2 \leq i \leq s$ and that $|B \cap C'| = \chi(C')$ for each $C' \in \C_{i-1}$. For any $\C_i$-class $C$ the set $B \cap C$ is the disjoint union of its subsets $B \cap C'$ over all $\C_{i-1}$-classes $C' \subseteq C$, that is to say, over all $C' \in \C_{i-1}(C)$. Then for each $C \in \C_i$,
    \[ 
    |B \cap C| = \sum_{C' \in \C_{i-1}(C)} |B \cap C'| = \sum_{C' \in \C_{i-1}(C)} \chi(C'), 
    \]
which by \eqref{array2} is equal to $\chi(C)$. Thus, by induction, $\chi(C) = |B \cap C|$ for all $i \in \{1, \ldots, s\}$ and all $C \in \C_i$. Hence $\chi = \chi_B$.
\end{proof}

From now on we usually omit the adjective `abstract' and  speak simply of `array functions'.
Next we define a group action on the set of array functions with respect to some $s$-chain of partitions, and use this to define an equivalence relation on array functions.

\begin{definition}\label{def:equivarray}
Let $G \leq \Sym(\PP)$ such that $G$ is contained in the stabiliser of an $s$-chain $\Ch$ of partitions as in \eqref{eq:chain}. Let $\chi : \bigcup_{i=1}^{s} \C_i \rightarrow \mathbb{Z}_{\geq 0}$ be an array function with respect to $\Ch$ and let $g\in G$, and define the function $\chi^g : \bigcup_{i=1}^{s} \C_i \rightarrow \mathbb{Z}_{\geq 0}$ by
    \[ 
    \chi^g(C) = \chi\big(C^{g^{-1}}\big) \quad \text{for all} \ C \in \bigcup_{i=1}^{s} \C_i. 
    \]
Define two array functions $\chi$ and $\chi'$ with respect to $\Ch$ to be \emph{$G$-equivalent} if $\chi^g = \chi'$ for some $g \in G$.
\end{definition}

\begin{lemma} \label{lem:array-equiv}
Let $s\geq2$, and suppose that $\Ch$ is an $s$-chain of partitions of $\PP$ as in \eqref{eq:chain}, with stabiliser $G = S_{e_1} \wr \ldots \wr S_{e_s}$.
\begin{enumerate}[{\rm (a)}]
    \item If $\chi : \bigcup_{i=1}^{s} \C_i \rightarrow \mathbb{Z}_{\geq 0}$ is an array function and $g\in G$, then also $\chi^g$ is an array function, and moreover, if $\chi=\chi_B$, where $B\subseteq\PP$, then $\chi^g=\chi_{B^g}$.
    \item The array functions $\chi_B$ and $\chi_{B'}$, where $B, B' \subseteq \PP$, are $G$-equivalent if and only if $B' = B^g$ for some $g \in G$. In particular, $G$-equivalence is an equivalence relation on the set of array functions with respect to $\Ch$.
\end{enumerate}
\end{lemma}

\begin{proof}
(a) Since, for each $i\in\{1,\dots,s-1\}$, $g$ induces a permutation of the classes of $\C_i$, condition \eqref{array1} follows immediately for $\chi^g$. Also, for $i\in\{2,\dots,s\}$ and $C\in\C_i$, $g^{-1}$ maps the set $\C_{i-1}(C)$ of $\C_{i-1}$-classes contained in $C$ to the set $\C_{i-1}(C^{g^{-1}})$ of $\C_{i-1}$-classes contained in $C^{g^{-1}}$. Hence, using condition \eqref{array2} for $\chi$, we have 
    \[
    \sum_{C' \in \C_{i-1}(C)} \chi^g(C') = \sum_{C' \in \C_{i-1}(C)} \chi((C')^{g^{-1}}) = \sum_{(C')^{g^{-1}} \in \C_{i-1}(C^{g^{-1}})} \chi((C')^{g^{-1}})  = \chi(C^{g^{-1}}) = \chi^g(C),
    \]
so \eqref{array2} holds also for $\chi^g$. Thus $\chi^g$ is an array function, by Definition~\ref{def:array}. Finally, suppose that $\chi=\chi_B$ for some $B\subseteq\PP$. Then, noting that $g^{-1}$ induces a permutation of $\bigcup_{i=1}^{s} \C_i$, for each $C\in \bigcup_{i=1}^{s} \C_i$, $\chi^g(C) = \chi(C^{g^{-1}}) = |B\cap C^{g^{-1}}| = |B^g\cap C|$, and hence $\chi^g = \chi_{B^g}$, proving part (a).

(b) Let $B, B' \subseteq \PP$. By part (a), $\chi_B^g=\chi_{B^g}$, and hence, if $B'=B^g$ then $\chi_B^g=\chi_{B'}$ and so $\chi_B$ and $\chi_{B'}$ are $G$-equivalent. Conversely suppose that $\chi_B$ and $\chi_{B'}$ are $G$-equivalent, so $\chi_B^g=\chi_{B'}$ for some $g\in G$, by Definition~\ref{def:equivarray}.
By part (a) this means that $\chi_{B^g}=\chi_{B'}$. For any subset $X\subseteq\PP$ and class $C\in\C_1$, let $X_C=X\cap C$. Then as we noted in the proof of Lemma~\ref{lem:array}, $B^g=\bigcup_{C\in\C_1}(B^g\cap C)$ and
$B'=\bigcup_{C\in\C_1}(B'\cap C)$, in each case a disjoint union. Also, for each $C\in\C_1$, 
    \[
    |B^g\cap C|=\chi_{B^g}(C)=\chi_{B'}(C)=|B'\cap C|.
    \]
As discussed in the paragraphs before Theorem~\ref{thm:BPRSchain}, the kernel of the action of $G$ on $\C_1$ is $G_{(\C_1)}\cong G_1^{|\C_1|}$, that is to say, this kernel is a direct product of $|\C_1|$ copies of $G_1=S_{e_1}$, with each copy acting naturally on one class of $\C_1$ and fixing all the other classes pointwise. For each $C\in\C_1$, since $|B^g\cap C|=|B'\cap C|$, 
there exists a permutation $h_C$ in the copy of $S_{e_1}$ acting on $C$ such that $h_C: B^g\cap C \to B'\cap C$. Let $h$ be the element of $G_{(\C_1)}$ which induces $h_C$ on $C$, for each $C\in\C_1$. Then $h$ maps $B^g$ to $B'$, and hence we have $gh\in G$ such that $B^{gh}=B'$ as required. This proves the first assertion of part (b). Finally, it is easy to see that $G$-equivalence is an equivalence relation on the set of array functions. The main issue is to use part (b) to prove transitivity of the relation (the other conditions are easier): if $\chi_B$ is $G$-equivalent to $\chi_{B'}$, and $\chi_{B'}$ is $G$-equivalent to $\chi_{B''}$, then by part (b) there exist $g, h\in G$ such that $B'=B^g$ and $B''=(B')^h$. This implies that $B''=B^{gh}$, and hence that $\chi_{B}$ is $G$-equivalent to $\chi_{B''}$.
\end{proof}

\section{Proofs of Theorems \ref{mainthm:2des} and \ref{thm:subdes}}

Let $s\geq2$, and suppose that the partitions  $\C_0, \ldots, \C_{s}$ of $\PP$ form an $s$-chain $\Ch$ as in \eqref{eq:chain} in Definition~\ref{def:chain}, with stabiliser $G = S_{e_1} \wr \ldots \wr S_{e_s}$. Consider a subset $B\subseteq\PP$ of size $k$, with array function $\chi_{B}$ with respect to $\Ch$, and let $\B = B^G$. It follows from Lemma \ref{lem:array-equiv} that $\B$ consists of all $k$-subsets $B'$ of $\PP$ such that the array function $\chi_{B'}$ with respect to $\Ch$ is $G$-equivalent to $\chi_B$. Recall from \eqref{cd} that the parameters $c_i$ and $d_i$ can be obtained from the numbers $e_i$, for $1 \leq i \leq s$. Thus $\D = (\PP,\B)$ is determined by the parameters $e_1, \ldots, e_s$, and an array function $\chi$ (up to $G$-equivalence) satisfying the conditions \eqref{array1} and \eqref{array2} of Definition \ref{def:array}. We denote $\D$ by $\D(e_1, \ldots, e_s; \,\chi)$, noting that this is the same notation used in \cite{CP93} when $s=2$.

In the case where $s=2$ the result  \cite[Proposition 2.2 (ii)]{CP93} gives necessary and sufficient conditions on $c$, $d$, and $\mathbf{x}$ in order for $\D(c,d;\mathbf{x})$ to be a $2$-design and to be a $3$-design; designs constructed in this manner are now called Cameron--Praeger designs. We generalise this result for $\D(e_1, \ldots, e_s; \,\chi)$ with $s \geq 2$ in Theorem \ref{mainthm:2des}, which we prove below. For any set $X$, the symbol $X^{\{2\}}$ denotes the set of all unordered pairs of elements of $X$.

Given an array function $\chi$, for brevity we will frequently use the notation
    \begin{equation} \label{x_C}
    x_C := \chi(C)
    \end{equation}
for any $C \in \C_i$, $1 \leq i \leq s$. 
If $B \subseteq \PP$ such that $\chi_B = \chi$, we can extend the notation \eqref{x_C} by defining $x_C$, for any $C = \{(\delta_j)_{j>0}\} \in \C_0$, to be
    \[ x_C = |B \cap C| = \begin{cases} 1 &\text{if } (\delta_j)_{j>0} \in B \\ 0 &\text{otherwise}. \end{cases} \]

\begin{proof}[Proof of Theorem \ref{mainthm:2des}]
Let $G, \PP$ be as above with $v=|\PP|=\prod_{i=1}^{s}e_i$, and let $B\subseteq\PP$ of size $k$. Let $\chi=\chi_B$ be the array function of $B$ with respect to $\Ch$ and recall the notation $x_C$ introduced above. Let $O_1, \ldots, O_s$ be as in \eqref{O_i}. By \cite[Proposition 1.3]{CP93}, the point-block incidence structure $\D = \big( \PP, B^G \big)$ is a $2$-design if and only if there is a constant $n$ such that $\big| B^{\{2\}} \cap O_i \big| = n\,|O_i|$ for all $i \in \{1, \ldots, s\}$. Assume such a $n$ exists. Note that $\sum_{i=1}^s |O_i| = \big| \PP^{\{2\}} \big| = \frac{1}{2} v(v-1)$, and $\sum_{i=1}^s \big| B^{\{2\}} \cap O_i \big| = \big| B^{\{2\}} \big| = \frac{1}{2} k(k-1)$.  Thus,  summing $n\,|O_i|$ over all $i$, we obtain
    \[ 
    n \cdot \frac{v(v-1)}{2} = \sum_{i=1}^s n \,|O_i| = \sum_{i=1}^s \big| B^{\{2\}} \cap O_i \big| = \frac{k(k-1)}{2}, 
    \]
so that $n = \frac{k(k-1)}{v(v-1)}$. Hence $\D$ is a $2$-design if and only if 
    \begin{equation} \label{ratio}
    \big| B^{\{2\}} \cap O_i \big| = \frac{k(k-1)}{v(v-1)} \cdot |O_i| \text{ for all } i \in \{1, \ldots, s\}.
    \end{equation}
Replacing $|O_i|$ with the formula in \eqref{|O_i|} we get
    \begin{equation} \label{ratio2}
    \frac{k(k-1)}{v(v-1)} \cdot |O_i| = \frac{k(k-1)}{v(v-1)} \cdot \frac{1}{2} \,v(e_i - 1) \prod_{j \leq i-1} e_j = \frac{1}{2} \cdot \frac{k(k-1)}{v-1} \,(e_i - 1) \prod_{j \leq i-1} e_j.
    \end{equation}
Recall from Section~\ref{s:innout} that the orbit $O_1$ consists of all $1$-inner pairs. Hence
    \[ 
    \big| B^{\{2\}} \cap O_1 \big|
    =  \sum_{C \in \C_1} \big| (B \cap C)^{\{2\}} \big|
    = \frac{1}{2} \sum_{C \in \C_1} |B \cap C| \cdot \big( |B \cap C| - 1\big)
    = \frac{1}{2} \sum_{C \in \C_1} x_C (x_C - 1). \]
To determine $\big| B^{\{2\}} \cap O_i \big|$ for $i \geq 2$, recall that point pairs $\{\alpha,\beta\}$ in $O_i$, that is,  $i$-strong inner pairs, must lie in the same $\C_i$-class and in different $\C_{i-1}$-classes, and observe that for any class $C' \in \C_i$, a pair $\{\alpha,\beta\}$ of points in $B \cap C'$ lies in the orbit $O_i$ exactly when $\alpha \in C$ for some $\C_{i-1}$-class $C \subseteq C'$ but $\beta \notin C$. With $\C_{i-1}(C')$ as in \eqref{C(X)} we have
    \begin{align*}
    \big| (B \cap C')^{\{2\}} \cap O_i \big|
    &= \frac{1}{2} \sum_{C \in \C_{i-1}(C')} |B \cap C| \cdot \big( |B \cap C'| - |B \cap C| \big) \\
    &= \frac{1}{2} \sum_{C \in \C_{i-1}(C')} x_{C} \big( x_{C'} - x_{C} \big)
    \end{align*}
and so
    \begin{align*}
    \big| B^{\{2\}} \cap O_i \big|
    &= \sum_{C' \in \C_i} \big| (B \cap C')^{\{2\}} \cap O_i \big|
    = \frac{1}{2} \sum_{C' \in \C_i} \ \sum_{C \in \C_{i-1}(C')} x_{C} \big( x_{C'} - x_{C} \big).
    \end{align*}
Note that the double summation runs over all pairs $(C',C)$ such that $C'\in\C_i$, $C\in\C_{i-1}$, and $C\subset C'$. For each $C\in\C_{i-1}$, the unique possibility for $C'$ is the class $C' = C^+$ 
defined in \eqref{C+}. Thus writing the double summation as a single summation over $C \in \C_{i-1}$, and writing $C^+$ instead of $C'$, we have
    \[
    \frac{1}{2} \sum_{C' \in \C_i} \ \sum_{C \in \C_{i-1}(C')} x_{C} \big( x_{C'} - x_{C} \big)
    = \frac{1}{2} \sum_{C \in \C_{i-1}} x_{C} \big( x_{C^+} - x_{C} \big).
    \]
Therefore, rewriting the conditions in \eqref{ratio} and \eqref{ratio2},  $\D$ is a $2$-design if and only if
    \[
    \frac{1}{2} \sum_{C \in \C_1} x_C (x_C - 1)
    = \frac{1}{2} \cdot \frac{k(k-1)}{v-1} \,(e_1 - 1)
    \]
and
    \[
    \frac{1}{2} \sum_{C \in \C_{i-1}} x_{C} \big( x_{C^+} - x_{C} \big)
    = \frac{1}{2} \cdot \frac{k(k-1)}{v-1} \,(e_i - 1) \prod_{j \leq i-1} e_j,\quad \text{ for all } i \in \{2, \ldots, s\}, 
    \]
which yield condition \eqref{2des-1} and conditions \eqref{2des} for all $i \in \{2, \ldots, s\}$.

Observe that the right hand side of the equalities \eqref{2des-1} and \eqref{2des} are even integers, since the expressions on the left hand side count the number of ordered pairs of distinct points that lie in the same $\C_1$-class for \eqref{2des-1}, and that lie in the same $\C_i$-class but different $\C_{i-1}$-classes for \eqref{2des}. Hence, for each $i \in \{1, \ldots, s\}$, the number $\frac{1}{2} \cdot \frac{k(k-1)}{v-1} \,(e_i - 1) \prod_{j \leq i-1} e_j = \frac{1}{v-1} \cdot \binom{k}{2} \,(e_i - 1) \prod_{j \leq i-1} e_j$ is an integer. Since $v-1 = \big(\prod_{j=0}^s e_j\big) - 1$ is coprime to $\prod_{j \leq i-1} e_j$, the number $v-1$ must divide $\binom{k}{2}(e_i - 1)$ for each $i \in \{1, \ldots, s\}$. Therefore $v-1$ divides $\binom{k}{2} \cdot \gcd(e_1-1, \ldots, e_s-1)$.

This proves Theorem~\ref{mainthm:2des}.
\end{proof}

We explore these conditions further.
Recall from Remark \ref{rem:s-cond} that conditions \eqref{2des-1} and \eqref{2des} consist of $s$ equations. We see in the next result that these $s$ equations are not independent.

\begin{lemma} \label{lem:2des-red}
If any $s-1$ equations out of the $s$ equations in \eqref{2des-1} and \eqref{2des} are satisfied, then the remaining equation is also satisfied.
\end{lemma}

\begin{proof}
Recall from the proof of Theorem \ref{mainthm:2des} that the condition \eqref{2des-1} is equivalent to \eqref{ratio} with $i=1$, and for any $i \in \{2, \ldots, s\}$ the condition \eqref{2des} is equivalent to \eqref{ratio}. Hence to prove the result it is enough to show that if \eqref{ratio} holds for all $i \in \{1, \ldots, s\} \setminus \{i'\}$ (for some $1 \leq i' \leq s$) then \eqref{ratio} holds for $i = i'$. Indeed, observe that
    \[ \big| B^{\{2\}} \cap O_{i'} \big| = \sum_{1 \leq i \leq s} \big| B^{\{2\}} \cap O_i \big| - \sum_{1 \leq i \leq s, \ i \neq i'} \big| B^{\{2\}} \cap O_i \big|. \]
Since $\sum_{i=1}^s \big| B^{\{2\}} \cap O_i \big| = \big| B^{\{2\}} \big| = \frac{1}{2} k(k-1)$, and since, by hypothesis, \eqref{ratio} holds for all $i \neq i'$, we can replace the terms of the sums on the right side of the equation above to obtain
    \begin{align*}
    \sum_{1 \leq i \leq s} \big| B^{\{2\}} \cap O_i \big| - \sum_{1 \leq i \leq s, \ i \neq i'} \big| B^{\{2\}} \cap O_i \big|
    &= \frac{k(k-1)}{2} - \sum_{1 \leq i \leq s, \ i \neq i'} \frac{k(k-1)}{v(v-1)} \cdot |O_i| \\
    &= \frac{k(k-1)}{2} - \frac{k(k-1)}{v(v-1)} \sum_{1 \leq i \leq s, \ i \neq i'} |O_i| \\
    &= \frac{k(k-1)}{2} - \frac{k(k-1)}{v(v-1)} \left( \left( \sum_{1 \leq i \leq s} |O_i| \right) - |O_{i'}| \right).
    \end{align*}
Now $\sum_{i=1}^s |O_i| = \big| \PP^{\{2\}} \big| = \frac{1}{2} v(v-1)$, so
    \[ \frac{k(k-1)}{v(v-1)} \left( \left( \sum_{1 \leq i \leq s} |O_i| \right) - |O_{i'}| \right) 
    = \frac{k(k-1)}{v(v-1)} \left( \frac{v(v-1)}{2} - |O_{i'}| \right)
    = \frac{k(k-1)}{2} - \frac{k(k-1)}{v(v-1)} \cdot |O_{i'}|. \]
Thus
    \[ \big| B^{\{2\}} \cap O_{i'} \big| = \sum_{1 \leq i \leq s} \big| B^{\{2\}} \cap O_i \big| - \sum_{1 \leq i \leq s, \ i \neq i'} \big| B^{\{2\}} \cap O_i \big| = \frac{k(k-1)}{v(v-1)} \cdot |O_{i'}|, \]
which shows that \eqref{ratio} holds for $i = i'$, as required. This proves the lemma.
\end{proof}

Corollary \ref{cor:2des-red} is an immediate consequence of Theorem \ref{mainthm:2des} and Lemma \ref{lem:2des-red}. Note that since $G = S_{e_1} \wr \ldots \wr S_{e_s}$, the structure $\D = \big(\PP, B^G\big)$ in Theorem \ref{mainthm:2des} is $\D(e_1, \ldots, e_s; \,\chi)$ for any array function $\chi$ that is $G$-equivalent to $\chi_B$.

\begin{corollary} \label{cor:2des-red}
Given integers $s \geq 2$ and $e_1, \ldots, e_s \geq 2$, and an array function $\chi : \bigcup_{i=1}^{s}  \C_i \rightarrow \mathbb{Z}_{\geq 0}$ as in Definition $\ref{def:array}$, the incidence structure $\D(e_1, \ldots, e_s; \,\chi)$ is a $2$-design if and only if $s-1$ out of the $s$ equations in \eqref{2des-1} and \eqref{2des} hold.
\end{corollary}

\begin{remark} \label{rem:2des-b}
Let $(\delta_j)_{j>i-1} \in \prod_{j>i-1} \Delta_j$. 
Recall from \eqref{partition} that $C_{(\delta_j)_{j>i-1}}$ is a class of $\C_{i-1}$, and from \eqref{C+} that $(C_{(\delta_j)_{j>i-1}})^+ = C_{(\delta_j)_{j>i}}$ is a class $C$ of $\C_i$. Using this labelling of $\C_i$, the parameter $x_C$ defined in \eqref{x_C} can be denoted by $x_C = x_{(\delta_j)_{j>i}}$. That is,
     \begin{equation} \label{array-specific} x_{(\delta_j)_{j>i}} := \big| B \cap C_{(\delta_j)_{j>i}}\big|, \end{equation}
so $x_{(\delta_j)_{j>i}} = \chi\big(C_{(\delta_j)_{j>i}}\big)$ for $i \geq 0$. Thus, using the notation \eqref{array-specific}, we can write the left hand side of \eqref{2des} as
    \[
    \sum_{C \in \C_{i-1}} x_C \left( x_{C^+} - x_C \right)
    = \sum_{(\delta_j)_{j>i-1} \in \prod_{j>i-1} \Delta_j} x_{(\delta_j)_{j>i-1}} \big( x_{(\delta_j)_{j>i}} - x_{(\delta_j)_{j>i-1}} \big).
    \]
Hence condition \eqref{2des} is equivalent to
    \begin{equation} \label{2des-b}
    \sum_{(\delta_j)_{j>i-1} \in \prod_{j>i-1} \Delta_j} x_{(\delta_j)_{j>i-1}} \left( x_{(\delta_j)_{j>i}} - x_{(\delta_j)_{j>i-1}} \right)
    = \frac{k(k-1)}{v-1} (e_i - 1) \prod_{j \leq i-1} e_j
    \end{equation}
    for $i\in\{2,\dots,s\}$.
\end{remark}

\begin{remark}
Assume that $s = 2$. By Corollary \ref{cor:2des-red}, the point-block structure $\D = \D(e_1, e_2; \,\chi)$ is a $2$-design if and only if condition \eqref{2des} holds for $i = 1$. It follows from \eqref{2des-1} that $\D$ is a $2$-design if and only if
    \[ 
    \sum_{\delta_2 \in \Delta_2} x_{(\delta_2)} \big( x_{(\delta_2)} - 1 \big) = \frac{k(k-1)}{v-1} (e_1 - 1). \]
Recall from Subsection \ref{subsec:chain} that $e_1 = c_1$ and that $|\Delta_2| = e_2 = d_1$; take $\Delta_2 = \{1, \ldots, d_1\}$. Thus $\D$ is a $2$-design if and only if
    \[ \sum_{\delta_2 = 1}^{d_1} x_{(\delta_2)} \big( x_{(\delta_2)} - 1 \big) = \frac{k(k-1)}{v-1} (c_1 - 1), \]
which is precisely  \cite[Proposition 2.2 (ii)]{CP93}.
\end{remark}

The proof of Theorem \ref{thm:subdes} uses the following result from \cite{CP93}.

\begin{proposition} \label{thm:strategy} \cite[Proposition 1.1]{CP93}
Let $\D = (\PP,\B)$ be a $t$-$(v,k,\lambda)$ design admitting a block-transitive (respectively, flag-transitive) group $G$. Let $H$ be a permutation group with $G \leq H \leq \Sym(\PP)$, and $\B^* = \big\{ B^H \ \big| \ B \in \B \big\}$. Then $\D^* = (\PP,\B^*)$ is a $t$-$(v,k,\lambda^*)$ design for some $\lambda^*$, admitting $H$ as a block-transitive (respectively, flag-transitive) group.
\end{proposition}

\begin{proof}[Proof of Theorem $\ref{thm:subdes}$]
It follows from Proposition \ref{thm:strategy} that there is a design $\D^* = (\PP,B^*)$ with the same point set $\PP$ as $\D$, with block set $\B^*$ containing the block set $\B$ of $\D$, such that $\D^*$ is $H$-block-transitive and $(H,s)$-chain-imprimitive, for $H = S_{e_1} \wr \ldots \wr S_{e_s}$. The design $\D^*$ has the same parameters $e_1, \ldots, e_s$, $v$, and $k$ as $\D$, and by Theorem \ref{mainthm:2des} these parameters must satisfy conditions \eqref{2des-1} and \eqref{2des}. This proves Theorem \ref{thm:subdes}. 
\end{proof}

\section{Constructions of chain-imprimitive designs} \label{sec:ex}

In this section we give in Construction \ref{con:s>2} an infinite family of  block-transitive, $s$-chain-imprimitive $2$-designs for any $s \geq 2$. Namely, for each $s$ our construction produces an infinite family of designs, with each design corresponding to an integer $p \geq 2$.

Before we give this construction, we first observe that examples of block-transitive, $s$-chain-imprimitive $2$-designs can be found among the Desarguesian projective planes. Indeed, suppose that $q$ is a prime power such that $q^2 + q + 1 = \prod_{i=1}^s e_i$ for some pairwise coprime integers $e_1, \ldots, e_s$, with $e_i \geq 2$ for each $i$. Let $\D = (\PP,\B)$ where $\PP$ and $\B$ are the point and block sets, respectively, of the projective plane ${\rm PG}(2,q)$. (We usually refer to the blocks of ${\rm PG}(2,q)$ as lines.) Then $\D$ is a $2$-$(v,k,1)$ design where $k = q + 1$ and $v = |\PP| = q^2 + q + 1$. Let $G$ be a subgroup of automorphisms of $\D$ generated by a Singer cycle. Then $G \cong \mathbb{Z}_{q^2 + q + 1}$ and acts transitively on blocks. Since the numbers $e_i$ are pairwise coprime, we have $G \cong \prod_{i=1}^s \mathbb{Z}_{e_i}$. The group $G$ leaves invariant the partitions $\C_1, \ldots, \C_s$, where $\C_j$ consists of the orbits of the subgroup $\prod_{i=1}^j \mathbb{Z}_{e_i}$ for each $j \in \{1, \ldots, s\}$. We have $\C_s = \{\PP\}$, and for $1 \leq j \leq s-1$ the partition $\C_j$ is nontrivial and is a proper refinement of $\C_{j+1}$. Thus these partitions form an $s$-chain and are preserved by $G$. 

Hence, in order to show that a block-transitive, $s$-chain-imprimitive projective plane exists for any integer $s$, it is sufficient to show that there does indeed exist a prime power $q$ such that $q^2 + q + 1$ can be written as a product of $s$ factors, each greater than $1$, which are pairwise coprime. We found that an exhaustive search over small values of $q$, using Magma \cite{magma}, did not readily produce examples of such factorisations $q^2+q+1=\prod_{i=1}^s e_i$ with large $s$, but we were able to spot a pattern which we developed into an infinite family of possibilities in Lemma \ref{lem:e_i}. 

\begin{lemma} \label{lem:e_i}
Let $q := p^{2^{s-1}}$ for positive integers $p > 1$ and $s \geq 3$. Let
    \[ e_1 := p^2 + p + 1 \]
and
    \[ e_i := p^{2^{i-1}} - p^{2^{i-2}} + 1 \quad \text{for $2 \leq i \leq s$}. \]
Then for any $j \in \{1, \ldots, s\}$:
    \begin{enumerate}[(i)]
    \item \label{e-prod} $\prod_{i=1}^j e_i = p^{2^j} + p^{2^{j-1}} + 1$ (in particular $\prod_{i=1}^s e_i = q^2 + q + 1$); and
    \item \label{e-coprime} $(e_i,e_j) = 1$ for any $i \neq j$.
    \end{enumerate}
\end{lemma}

\begin{proof}
Clearly $e_1 = p^2 + p + 1 = p^{2^1} + p^{2^0} + 1$. Arguing by induction, for $j \geq 2$ we have
    \begin{align*}
    \prod_{i=1}^j e_i = \left(\prod_{i=1}^{j-1} e_i\right) e_j
    &= \big( p^{2^{j-1}} + p^{2^{j-2}} + 1 \big)\big( p^{2^{j-1}} - p^{2^{j-2}} + 1 \big) \\
    &= \frac{\big( p^{2^{j-2}} \big)^3 - 1}{p^{2^{j-2}} - 1} \cdot \frac{\big( p^{2^{j-2}} \big)^3 + 1}{p^{2^{j-2}} + 1} \\
    &= \frac{\big( p^{2^{j-2}} \big)^6 - 1}{\big(p^{2^{j-2}}\big)^2 - 1} \\
    &= \big( p^{2^{j-2}} \big)^4 + \big( p^{2^{j-2}} \big)^2 + 1 \\
    &= p^{2^j} + p^{2^{j-1}} + 1.
    \end{align*}
This proves part \eqref{e-prod}.

We now prove part \eqref{e-coprime}. Assume first that either $i$ or $j$ is $1$; without loss of generality suppose that $i = 1$. Let $D = (e_1, e_j) = \big( p^2 + p + 1, \ p^{2^{j-1}} - p^{2^{j-2}} + 1 \big)$. Then $D$ divides $\big( p^3 - 1, \ p^{3 \cdot 2^{j-2}} + 1 \big)$. Now
    \[ \big( p^3 - 1, \ p^{3 \cdot 2^{j-2}} + 1 \big) = \big( p^3 - 1, \ \big( p^{3 \cdot 2^{j-2}} - 1 \big) + 2 \big) = (p^3 - 1, \ 2).  \]
Hence $D = 1$ or $2$; since $p^2 + p + 1$ is odd for any $p$, we must have $D = 1$.

The case where $i,j \geq 2$ is proved by a similar argument. Without loss of generality suppose that $i < j$. Then $D := (e_i, e_j) = \big( p^{2^{i-1}} - p^{2^{i-2}} + 1, \ p^{2^{j-1}} - p^{2^{j-2}} + 1 \big)$ divides $\big( p^{3 \cdot 2^{i-2}} + 1, \ p^{3 \cdot 2^{j-2}} + 1 \big) = \big( p^{3 \cdot 2^{i-2}} + 1,  \ \big( p^{3 \cdot 2^{j-2}} - 1 \big) + 2 \big)$. Now $2^{j-2} = 2^{i-2} \cdot 2^{j-i}$ and $p^{3 \cdot 2^{j-2}} - 1 = \big(p^{3 \cdot 2^{i-2}}\big)^{2^{j-i}} - 1$. Since $2^{j-i}$ is even, $\big(p^{3 \cdot 2^{i-2}}\big)^{2^{j-i}} - 1$ is divisible by $p^{3 \cdot 2^{i-2}} + 1$. Thus
    \[ \big( p^{3 \cdot 2^{i-2}} + 1,  \ \big( p^{3 \cdot 2^{j-2}} - 1 \big) + 2 \big) = \big( p^{3 \cdot 2^{i-2}} + 1, \ 2 \big) \]
and again $D = 1$ or $2$. Since $p^{2^{i-1}} - p^{2^{i-2}} + 1$ is odd for any $p$ and $i$, it follows that $D = 1$. This completes the proof of part \eqref{e-coprime}.
\end{proof}

It follows from Lemma \ref{lem:e_i} that for a given $s \geq 2$, taking $q = p^{2^{s-1}}$ with $p$ a prime power gives us a block-transitive, $s$-chain-imprimitive design from the projective plane ${\rm PG}(2,q)$. We used these examples, in particular ${\rm PG}(2,16)$, as a starting point for our construction, but we note that our construction works for all integers $p$, not only prime powers. The general strategy that we use was inspired by \cite[Proposition 1.1]{CP93} (see Theorem \ref{thm:strategy}). By Theorem \ref{thm:strategy}, taking any block $B$ in a $G$-block-transitive, $s$-chain-imprimitive projective plane, where $G$ is a Singer cycle, and replacing the block set with $B^H$, where $G< H = S_{e_1} \wr \ldots \wr S_{e_s}$, yields another block-transitive, $s$-chain-imprimitive $2$-design. This is precisely the strategy that we use to obtain an explicit construction. Our first attempt at a general construction began with  the projective plane $\PG(2,16)$, where we had $1+16+256 = 7\times 3\times 13$, so $s = 3$. We took a block $B$ and a group $G$ generated by a Singer cycle, and used \textsc{Magma} \cite{magma} to determine the `array' of $B$ with respect to the $G$-invariant partitions $\C_1$ and $\C_2$, where $\C_1, \C_2$ were the sets of orbits of the subgroups $\mathbb{Z}_{7}$ and $\mathbb{Z}_{7} \times \mathbb{Z}_{3}$, respectively; the parameters are as in Lemma \ref{lem:e_i} with $p=2$, $s=3$, $e_1=7$, $e_2=3$ and $e_3=13$. Our construction makes use of the labelling of points introduced in Section~\ref{subsec:chain}, and in order to build some intuition about the construction it is helpful to visualise the case $s=3$. Thus we describe this case first in Section~\ref{sub:ex}. In these descriptions $\mathbb{Z}_{e}^*$ denotes the set of non-zero elements of $\mathbb{Z}_{e}$.

\subsection{Examples with $s=2$} These examples are precisely the designs from Construction~\ref{con:s>2} with $s=2$. In the notation of \cite{CP93} these are the designs $\D\big(e_1,e_2; \, (p+1, \, 1^{e_2-1})\big)$.

\begin{construction} \label{con:s=2}
For any integer $p \geq 2$ let $q = p^2$, $e_1 = p^2 + p + 1$, and $e_2 = p^2 - p + 1$. Let $\PP = \mathbb{Z}_{e_1} \times \mathbb{Z}_{e_2}$ and $G = S_{e_1} \wr S_{e_2} < \Sym(\PP)$, and define the incidence structure $\D^2 = (\PP,\B)$ with $\B = B_2^G$ where
    \begin{equation} \label{B:s=2}
    B_2:= \{ (\delta_1, 0) \ | \ 0 \leq \delta_1 \leq p \} \cup \{ (0, \delta_2) \ | \ \delta_2 \in \mathbb{Z}_{e_2}^* \}.
    \end{equation}
The group $G$ stabilises the partition $\C$ consisting of $e_2$ classes $C_{(\delta_2)} = \{ (\varepsilon_1,\delta_2) \in \PP \ | \ \varepsilon_1 \in \mathbb{Z}_{e_1} \}$, with each class of size $e_1$. 
The block $B_2$ consists of $p+1$ points from the class $C_{(0)}$, and one point from each of the remaining $e_2-1$ classes. It will be shown in Theorem \ref{thm:s>1} that $\D^2$ is a $G$-block-transitive, $(G,2)$-chain-imprimitive $2$-design with $v = e_1e_2 = p^4 + p^2 + 1$ and $k = p^2 + 1$. Hence $k(k-1) = v-1$, and the parameter $\lambda = \frac{bk(k-1)}{v(v-1)} = \frac{b}{v}$ where $b := |\B|$. By Lemma \ref{lem:array-equiv} the block set $\B$ consists of all $k$-subsets $B$ of $\PP$ whose array function $\chi_{B}$ is $G$-equivalent to $\chi_{B_2}$, that is, $B$ consists of $p+1$ points from one class of $\C$, and one point from each of the remaining $e_2-1$ classes. Thus
    \begin{align*}
    b
    &= (\#\text{$(p+1)$-subsets in a class}) \cdot |\C| \cdot (\#\text{points in a class})^{|\C|-1} \\
    &= \binom{e_1}{p+1} \cdot e_2 \cdot e_1^{e_2-1}
    = \binom{e_1-1}{p} \cdot \frac{e_1}{p+1} \cdot e_2 \cdot e_1^{e_2-1}
    = \binom{e_1-1}{p} \cdot \frac{v}{p+1} \cdot e_1^{e_2-1}
    \end{align*}
and therefore $\lambda = \binom{e_1-1}{p} \cdot \frac{1}{p+1} \cdot e_1^{e_2-1}$. For future use it will be convenient to express $b$ and $\lambda$ using the notation $c_i$ defined in \eqref{cd}. In particular $e_1 = c_1$ and $v=c_2$, so we may write:
    \[
   b=\binom{e_1-1}{p} \cdot \frac{c_2}{p+1} \cdot c_1^{e_2-1} \quad\text{and} \quad \lambda = \binom{e_1-1}{p} \cdot \frac{1}{p+1} \cdot c_1^{e_2-1}.
    \]

\end{construction}

\subsection{Examples with $s=3$}\label{sub:ex} These examples are precisely the designs from Construction~\ref{con:s>2} with $s=3$. For this small value of $s$ it is easy to give concrete representations of the base block (Figure~\ref{fig:s=3}) and tabular information about the parameters (Table~\ref{tab:s=3}). We separate out this case from Construction \ref{con:s>2} because it contains all the components of that construction, and understanding this case will help the reader to visualise the construction for larger values of $s$.

\begin{construction} \label{con:s=3}
For any integer $p \geq 2$, let $q = p^4$, $e_1 = p^2 + p + 1$, $e_2 = p^2 - p + 1$, and $e_3 = p^4 - p^2 + 1$. Let $\PP = \mathbb{Z}_{e_1} \times \mathbb{Z}_{e_2} \times \mathbb{Z}_{e_3}$ and $G = S_{e_1} \wr S_{e_2}\wr S_{e_3} < {\Sym}(\PP)$, and define the point-block incidence structure $\D^3 = (\PP,\B)$ with $\B = B_3^G$, where
    \begin{equation} \label{B:s=3}
    B_3 := \{ (\delta_1, 0, 0) \ | \ 0 \leq \delta_1 \leq p \} \cup \{ (0, \delta_2, 0) \ | \ \delta_2 \in \mathbb{Z}_{e_2}^* \} \cup \{ (0, 0, \delta_3) \ | \ \delta_3 \in \mathbb{Z}_{e_3}^* \}.
    \end{equation}

The group $G$ is the stabiliser of the $3$-chain of partitions $\PP = \C_0 \prec \C_1 \prec \C_2 \prec \C_3 = \{\PP\}$, where
    \begin{itemize}
    \item the $\C_1$-classes are the sets $C_{(\delta_2,\delta_3)} = \{ (\varepsilon_1,\delta_2,\delta_3) \in \PP \ | \ \varepsilon_1 \in \mathbb{Z}_{e_1} \}$ for each $(\delta_2,\delta_3) \in \mathbb{Z}_{e_2} \times \mathbb{Z}_{e_3}$;
    \item the $\C_2$-classes are the sets $C_{(\delta_3)} = \{ (\varepsilon_1,\varepsilon_2,\delta_3) \ | \ \varepsilon_1 \in \mathbb{Z}_{e_1}, \varepsilon_2 \in \mathbb{Z}_{e_2} \}$ for each $\delta_3 \in \mathbb{Z}_{e_3}$.
    \end{itemize}
Figure \ref{fig:s=3} illustrates the distribution of the points of $B_3$ among the $\C_1$- and $\C_2$-classes; the numbers indicate the number of points in $B_3$ that belong to that class of $\C_1$ or $\C_2$. The values of the parameters $x_{(\delta_2,\delta_3)} = |B_3 \cap C_{(\delta_2,\delta_3)}|$ and $x_{(\delta_3)} = |B_3 \cap C_{(\delta_3)}|$, for $C_{(\delta_2,\delta_3)} \in \C_1(B_3)$ and $C_{(\delta_3)} \in \C_2(B_3)$ (with $\C_1(B_3)$ and $\C_2(B_3)$ as defined in \eqref{C(X)})  are listed in the second and third columns, respectively, of Table \ref{tab:s=3}. We can deduce these numbers easily by referring to Figure \ref{fig:s=3}. For $(\delta_2,\delta_3) \in \mathbb{Z}^*_{e_2} \times \mathbb{Z}^*_{e_3}$ the parameters $x_C = x_{(\delta_2,\delta_3)}$ have value $0$, and thus the terms $x_C(x_C - 1)$ and $x_C(x_{C^+} - x_C)$ in the sums \eqref{2des-1} and \eqref{2des} are both equal to $0$; for this reason, we omit these values from the table. Note that $x_{(0)}$ is the number of points in the first column of Figure \ref{fig:s=3}, so $x_{(0)} = (p + 1) + (e_2 - 1) = p^2 + 1$. In our analysis of the general construction (Theorem \ref{thm:s>1}) we will show that $\D^3$ is a $G$-block-transitive and $(G,3)$-chain-imprimitive $2$-design.

\begin{table}[ht]
    \centering
    \begin{tabular}{lcccc}
    \hline
    $(\delta_2,\delta_3)$ & $x_{(\delta_2,\delta_3)}$ & $x_{(\delta_3)}$ & $x_{(\delta_2,\delta_3)} - 1$ & $x_{(\delta_3)}-x_{(\delta_2,\delta_3)} $ \\
    \hline\hline
    $(0,0)$ & $p+1$ & $p^2 + 1$ & $p$ & $p^2 - p$ \\
    $(\delta_2,0)$, $\delta_2 \in \mathbb{Z}^*_{e_2}$ & $1$ & $p^2 + 1$ & $0$ & $p^2$ \\
    $(0,\delta_3)$, $\delta_3 \in \mathbb{Z}^*_{e_3}$ & $1$ & $1$ & $0$ & $0$ \\
    \hline
    \end{tabular}
    \caption{Values of $x_C$ for $C \in \C_1(B_3)$ and $C \in \C_2(B_3)$ in Construction \ref{con:s>2} with $s=3$}
    \label{tab:s=3}
\end{table}

\def\chainblocktwo{
    \multido{\n=-1+1.0}{3}{\psline(0,\n)(1.5,\n)}
    \pspolygon(0,-2)(0,2)(6,2)(6,-2)
    \multido{\n=1.5+1.5}{3}{\psline(\n,-2)(\n,2)}
    \rput[b](0.75,2.25){$0$} \rput[b](2.25,2.25){$1$} \rput[b](3.75,2.25){$\ldots$} \rput[r](-0.25,1.5){$0$} \rput[r](-0.25,0.5){$1$} \rput[r](-0.25,-0.5){$\vdots$} }
\def\onescolchain{\rput(0.75,0.5){$1$} \rput(0.75,-0.5){$\vdots$} \rput(0.75,-1.5){$1$}}
\def\onesrowchain{\rput(2.25,1.5){$1$} \rput(3.75,1.5){$\ldots$} \rput(5.25,1.5){$1$}}

\begin{figure}[ht]
    \centering
    \begin{pspicture}(0,-2.5)(6.25,3.5)
    \chainblocktwo {\onescolchain}{\onesrowchain}
    \rput[b](5.25,2.25){$e_3-1$} \rput[bl](6.75,2.25){$\longleftarrow$ $\delta_3$-labels}
    \rput[r](-0.25,-1.5){$e_2-1$} \rput[br](-0.25,2.25){\parbox[c]{1.5cm}{\flushright$\delta_2$-labels \\ $\downarrow$}}
    \rput(0.75,1.5){$p+1$}
	\pnode(1.5,2.75){a} \pnode(3,2.75){b} \ncbar[nodesep=2pt,angle=90,armA=5pt,armB=5pt]{a}{b} \naput{$\C_2$-class $C_{(\delta_3)}$ (column)}
    \pnode(-0.75,1){c} \pnode(-0.75,0){d} \ncbar[nodesep=2pt,angle=180,armA=5pt,armB=5pt]{c}{d} \nbput{\parbox[c]{2cm}{\flushright$\C_1$-class \\ $C_{(\delta_2,\delta_3)}$}}
	\end{pspicture}
	\caption{Distribution of points of $B_3$ for Construction \ref{con:s>2} with $s=3$}
	\label{fig:s=3}
\end{figure}

Since $v = e_1e_2e_3 = p^8 + p^4 + 1$ and $k = p^4 + 1$, we again have $k(k-1) = v-1$ so that $\lambda = \frac{b}{v}$ where $b := |\B|$. By Lemma \ref{lem:array-equiv}(b) the block set $\B = B_3^G$ consists of all $k$-element subsets $B$ of $\PP$ whose array function $\chi_{B}$ is $G$-equivalent to $\chi_{B_3}$. In order to find $b$, observe that the intersection of the generating block $B_3$ in \eqref{B:s=3} and the $\C_2$-class $C_{(0)}$ (the first column in Figure \ref{fig:s=3}) is the set $B_2 \times \{0\}$, where $B_2$ is the generating block of the design in Example \ref{con:s=2} given in \eqref{B:s=2}. For any subset $B'_2 \in B_2^{S_{e_1} \wr S_{e_2}}$, the set $(B'_2 \times \{0\}) \cup \big\{ (0,0,\delta_3) \ | \ \delta_3 \in \mathbb{Z}^*_{e_3} \big\}$ is also a block in $\B$. Let $b_2 := \big|B_2^{S_{e_1} \wr S_{e_2}}\big|$; then
    \[
    b
    = b_2 \cdot |\C_2| \cdot (\#\text{points in a $\C_2$-class})^{e_3-1}
    = b_2 \cdot e_3 \cdot c_2^{e_3-1}.
    \]
We found in Example \ref{con:s=2} that $b_2 = \binom{e_1-1}{p} \cdot \frac{e_1e_2}{p+1} \cdot c_1^{e_2-1}$; substituting this into the above gives
    \[
    b = \binom{e_1-1}{p} \cdot \frac{v}{p+1} \cdot c_1^{e_2-1} \cdot c_2^{e_3-1}
    \]
and therefore $\lambda = \binom{e_1-1}{p} \cdot \frac{1}{p+1} \cdot c_1^{e_2-1} \cdot c_2^{e_3-1}$.
\end{construction}

\subsection{The general construction} \label{subsec:s>2} We give here the general design construction which has two integer parameters, namely $p\geq2$ and $s\geq 2$.

\begin{construction} \label{con:s>2}
Let $p$ be an integer (not necessarily prime) with $p \geq 2$, and let $q = p^{2^{s-1}}$, $e_1 = p^2 + p + 1$, and $e_i = p^{2^{i-1}} - p^{2^{i-2}} + 1$ for $2 \leq i \leq s$. Let $\PP = \prod_{i=1}^s \mathbb{Z}_{e_i}$ and $G = S_{e_1} \wr \ldots \wr S_{e_s}<{\rm Sym}(\PP)$, and define $\D^s = (\PP,\B)$ with $\B = B_s^G$, where 
    \begin{equation} \label{B:s>2}
    B_s := \big\{ (\delta_1,0,\ldots,0) \in \PP \ \big| \ 0 \leq \delta_1 \leq p \big\} \cup \left( \bigcup_{\ell=2}^s \big\{ (\underbrace{0,\ldots,0}_{\ell-1}, \delta_\ell, 0,\ldots,0) \in \PP \ \big| \ \delta_\ell \in \mathbb{Z}^*_{e_\ell} \big\} \right).
    \end{equation}
\end{construction}

The group $G$ is the stabiliser of the $s$-chain of partitions $\Ch$ defined in \eqref{eq:chain} such that, for $1 \leq i \leq s-1$, the $\C_i$-classes are the sets
    \[ C_{(\delta_j)_{j>i}} = \big\{ (\varepsilon_1, \ldots, \varepsilon_i, \delta_{i+1}, \ldots, \delta_s) \in \PP \ \big| \ \varepsilon_j \in \mathbb{Z}_{e_j} \ \forall\, 1 \leq j \leq i \big\} 
    \]
for each $(\delta_j)_{j>i} = (\delta_{i+1}, \ldots, \delta_s) \in \prod_{j>i} \mathbb{Z}_{e_j}$, as in \eqref{class}, that is to say, $(\varepsilon_j)_{j=1}^s \in C_{(\delta_j)_{j>i}}$ if and only if $\varepsilon_j = \delta_j$ for all $j >i$. The distribution of the points of $B_s$ among the $\C_i$-classes is as follows. Let $C = C_{(\delta_j)_{j>i}}$. If $(\delta_j)_{j>i} = (0)_{j>i}$, then
    \[ 
    B_s \cap C = \big\{ (\delta_1,0,\ldots,0) \in \PP \ \big| \ 0 \leq \delta_1 \leq p \big\} \cup \left( \bigcup_{\ell=2}^i \big\{ (\underbrace{0,\ldots,0}_{\ell-1}, \delta_\ell, 0,\ldots,0) \in \PP \ \big| \ \delta_\ell \in \mathbb{Z}^*_{e_\ell} \big\} \right) 
    \]
which has size
    \begin{equation} \label{x_0}
    x_{(0)_{j>i}} = (p+1) + \sum_{\ell=2}^i (e_\ell - 1) = (p+1) + \sum_{\ell=2}^i \big( p^{2^{\ell-1}} - p^{2^{\ell-2}} \big) = p^{2^{i-1}} + 1.
    \end{equation}
On the other hand, for each $i' > i$ and $\delta_{i'} \in \mathbb{Z}_{e_{i'}}$, if $(\delta_j)_{j>i} = (\underbrace{0, \ldots, 0}_{i'-i-1}, \delta_{i'}, 0, \ldots, 0)$, then
    \[ B_s \cap C = \big\{ (\underbrace{0, \ldots, 0}_{i'-1}, \delta_{i'}, 0, \ldots, 0) \big\} \]
and
    \begin{equation} \label{x_delta}
    x_{({\scriptsize\underbrace{0, \ldots, 0}_{i'-i-1}}, \delta_{i'}, 0, \ldots, 0)} = 1.
    \end{equation}
For all other $\C_i$-classes the label $(\delta_j)_{j>i}$ has at least two nonzero entries, so $B_s \cap C = \varnothing$ and $x_{(\delta_j)_{j>i}} = 0$. The second and third columns of Table \ref{tab:s>2} lists the values of $x_C$ for $C \in \C_{i-1}(B_s) \cup \C_{i}(B_s)$; the values of $x_{(\delta)_{j>i-1}}$ are obtained by replacing $i$ by $i-1$ in \eqref{x_0} and \eqref{x_delta}. The distribution of the points in the $\C_{i+1}$-class $C_{(0)_{j>i+1}}$, for $1 \leq i \leq s-1$, are illustrated in Figure \ref{fig:s>2}, showing the $\C_{i-1}$ and $\C_i$-classes. As in Figure \ref{fig:s=3} the number in each cell indicates the number of points of $B_s$ that belong in that partition class. For $(\delta_j)_{j>i-1}$ with more than one nonzero coordinate, the parameters $x_C = x_{(\delta_j)_{j>i-1}}$ have value $0$, and thus the terms $x_C(x_C - 1)$ and $x_C(x_{C^+} - x_C)$ in the sums \eqref{2des-1} and \eqref{2des} are both equal to $0$; hence we omit these values from Table \ref{tab:s>2}.

\begin{table}[!ht]
\centering
    \begin{tabular}{lcccc}
    \hline
    $(\delta_j)_{j>i-1}$ ($2 \leq i \leq s$) & $x_{(\delta_j)_{j>i-1}}$ & $x_{(\delta_j)_{j>i}}$ & $x_{(\delta_j)_{j>i-1}} - 1$ & $x_{(\delta_j)_{j>i}} - x_{(\delta_j)_{j>i-1}}$ \\
    \hline\hline
    $(\underbrace{0,\ldots,0}_{s-i+1})$ & $p^{2^{i-2}} + 1$ & $p^{2^{i-1}} + 1$ & $p^{2^{i-2}}$ & $p^{2^{i-1}} - p^{2^{i-2}}$ \\
    $(\delta_i,\underbrace{0,\ldots,0}_{s-i})$, $\delta_i \in \mathbb{Z}^*_{e_i}$ & $1$ & $p^{2^{i-1}} + 1$ & $0$ & $p^{2^{i-1}}$ \\
    $(\underbrace{0,\ldots,0}_{i'-i},\delta_{i'},\underbrace{0,\ldots,0}_{s-i'})$, \ $i+1 \leq i' \leq s$ & $1$ & $1$ & $0$ & $0$ \\
    \hline
    \end{tabular}
    \caption{Values of $x_C$ for $C \in \C_{i-1}(B_s)$ and $C \in \C_i(B_s)$ in Construction \ref{con:s>2}}
    \label{tab:s>2}
\end{table}

\def\chainblocktwo{
    \multido{\n=-1+1.0}{3}{\psline(0,\n)(1.5,\n)}
    \pspolygon(0,-2)(0,2)(6,2)(6,-2)
    \multido{\n=1.5+1.5}{3}{\psline(\n,-2)(\n,2)}
    \rput[b](0.75,2.25){$0$} \rput[b](2.25,2.25){$1$} \rput[b](3.75,2.25){$\ldots$} \rput[r](-0.25,1.5){$0$} \rput[r](-0.25,0.5){$1$} \rput[r](-0.25,-0.5){$\vdots$} }
\def\onescolchain{\rput(0.75,0.5){$1$} \rput(0.75,-0.5){$\vdots$} \rput(0.75,-1.5){$1$}}
\def\onesrowchain{\rput(2.25,1.5){$1$} \rput(3.75,1.5){$\ldots$} \rput(5.25,1.5){$1$}}

\begin{figure}[ht]
    \centering
    \begin{pspicture}(-1,-2.5)(7,3.5)
    \chainblocktwo {\onescolchain}{\onesrowchain}
    \rput[b](5.25,2.25){$e_{i+1}-1$} \rput[bl](6.75,2.25){$\longleftarrow$ $\delta_{i+1}$-labels}
    \rput[r](-0.25,-1.5){$e_{i}-1$} \rput[br](-0.25,2.25){\parbox[c]{2cm}{\flushright$\delta_{i}$-labels \\ $\downarrow$}}
    \rput(0.75,1.5){$p^{2^{i-2}}+1$}
	\pnode(1.5,2.75){a} \pnode(3,2.75){b} \ncbar[nodesep=2pt,angle=90,armA=5pt,armB=5pt]{a}{b} \naput{$\C_{i}$-class $C_{(\delta_j)_{j>i}}$ (column)}
    \pnode(-0.75,2){c} \pnode(-0.75,1){d} \ncbar[nodesep=2pt,angle=180,armA=5pt,armB=5pt]{c}{d} \nbput{\parbox[t]{2cm}{\flushright$\C_{i-1}$-class \\ $C_{(0)_{j>i-1}}$}}
    \end{pspicture}
    \caption{Distribution of points of $B_s$ inside the $\C_{i+1}$-class $C_{(0)_{j>i+1}}$, for the case $s \geq 3$}
    \label{fig:s>2}
\end{figure}

\begin{lemma} \label{lem:b}
Let $\D^s, p, q, s, G$ be as in Construction $\ref{con:s>2}$. Then
    \begin{equation} \label{eq:b}
    |\B| = \binom{e_1-1}{p} \cdot \frac{v}{p+1} \cdot \prod_{i=2}^s (c_{i-1})^{e_i-1}.
    \end{equation}
\end{lemma}

\begin{proof}
Our proof is by induction. Let $b_s = |\B|$. Recall from Example \ref{con:s=2} that
    \[
    b_2 = \binom{e_1-1}{p} \cdot \frac{c_2}{p+1} \cdot c_1^{e_2-1},
    \]
where $c_2 = e_1e_2 = |\mathbb{Z}_{e_1} \times \mathbb{Z}_{e_2}|$, so the formula \eqref{eq:b} holds for $s=2$. Let $s \geq 3$ and assume that \eqref{eq:b} holds for $\ell \in \{2, \ldots, s-1\}$. Observe that the intersection of the block $B_s$ in \eqref{B:s>2} and the $\C_{s-1}$-class $C_{(0)}$ is $B_{s-1} \times \{0\}$, where $B_{s-1}$ is the generating block of the design $\D^{s-1}$ in Construction \ref{con:s>2}. Since, by Lemma \ref{lem:array-equiv}, the block set $\B$ consists of all $k$-subsets of $\PP$ whose array function is equivalent to $\chi_{B_s}$, it follows that for any subset $B'_{s-1} \in B_{s-1}^{S_{e_1} \wr \ldots \wr S_{e_{s-1}}}$, the set $(B'_{s-1} \times \{0\}) \cup \big\{ (0,\ldots,0,\delta_s) \in \PP \ \vline \ \delta_s \in \mathbb{Z}_{e_s} \big\}$ is a block in $\B$. Then
    \[
    b_s
    = b_{s-1} \cdot |\C_{s-1}| \cdot (\#\text{points in a $\C_{s-1}$-class})^{e_s-1}
    = b_{s-1} \cdot e_s \cdot (c_{s-1})^{e_s-1}.
    \]
By the induction hypothesis, $b_{s-1} = \binom{e_1-1}{p} \cdot \frac{c_{s-1}}{p+1} \cdot \prod_{i=2}^{s-1} (c_{i-1})^{e_i-1}$. Substituting this into the previous equation, and recalling that $c_{s-1}e_s = c_s$, we obtain
    \[
    b_s
    = \binom{e_1-1}{p} \cdot \frac{c_{s-1}}{p+1} \cdot \prod_{i=2}^{s-1} (c_{i-1})^{e_i-1} \cdot e_s \cdot (c_{s-1})^{e_s-1}
    = \binom{e_1-1}{p} \cdot \frac{c_s}{p+1} \cdot \prod_{i=2}^s (c_{i-1})^{e_i-1}.
    \]
This completes the proof.
\end{proof}

\begin{theorem} \label{thm:s>1}
Let $\D^s, q, s, G$ be as in Construction $\ref{con:s>2}$. Then $\D^s$ is a $2$-design with $v = q^2 + q + 1$, $k = q+1$, and $\lambda = \binom{e_1-1}{p} \cdot \frac{1}{p+1} \cdot \prod_{i=2}^s (c_{i-1})^{e_i-1}$, $\D^s$ is not a $3$-design, and $G \leq \Aut(\D^s)$ is block-transitive and $s$-chain-imprimitive, but not flag-transitive.
\end{theorem}

\begin{proof}
 It follows from Lemma \ref{lem:e_i}\eqref{e-prod} that
    \[ v = |\PP| = \prod_{i=1}^s e_i = p^{2^s} + p^{2^{s-1}} + 1 = q^2 + q + 1, \]
and from the definition of $B_s$ in \eqref{B:s>2} that
    \[ 
    k = |B_s|  = (p+1) + \sum_{\ell=2}^s (e_\ell - 1) = (p+1) + \sum_{\ell=2}^s (p^{2^{\ell-1}} - p^{2^{\ell-2}})
    = p^{2^{s-1}} + 1 = q + 1. \]
By Construction  $\ref{con:s>2}$, it is clear that $G\leq \Aut(\D^s)$, and that $G$ is block-transitive and $s$-chain-imprimitive. However $G$ is not flag-transitive, since if it were then the block stabiliser $G_{B_s}$ would be transitive on $B_s$, whereas $G_{B_s}$ fixes setwise $B_s\cap C$, for $C=C_{(0)_{j>1}}$, the unique class of $\C_1$ containing more than one point of $B_s$. 

Next we prove that $\D^s$ is not a $3$-design. First we make an observation. Consider the following $(p+e_2)$-element subset of $\PP$:
$$
    Z := \big\{ (\delta_1,0,\ldots,0) \in \PP \ \big| \ 0 \leq \delta_1 \leq p \big\} \cup \big\{ (0, \delta_2, 0,\ldots,0) \big| \ \delta_2 \in \mathbb{Z}^*_{e_2} \big\}.
$$ 
Then $Z\subset B_s$, and we note that each pair of points from $Z$ is either a $1$-strong inner pair, or the pair contains a least one point of the form $(0, \delta_2, 0,\ldots,0)$ and is a $2$-strong inner pair (as defined in Section~\ref{s:innout}). Let $r'$ be the number of blocks of $\D^s$ containing $Z$. By the transitivity of $G$, if $Z'$ is a $(p+e_2)$-subset of $\PP$ which consists of a $(p+1)$-subset of some $\C_1$-class $C$, together with one point from each of the $\C_1$-classes contained in $C^+\setminus C$, where $C^+\in\C_2$ as in \eqref{C+}, then $Z'$ is contained in exactly $r'$ blocks of $\D^s$.   We now consider two $3$-subsets $X, Y$ of points: $X$ is contained in a $\C_1$-class $C$, while $Y=\{y_1,y_2,y_3\}$ where $\{y_1, y_2\}$ is contained in a $\C_1$-class $C$ and $y_3\in C^+\setminus C$. Each pair from $X$ is a $1$-strong inner pair, while $Y$ contains one $1$-strong inner pair and two $2$-strong inner pairs. The number $\lambda_X$ of blocks containing $X$ is equal to the number $\binom{e_1-3}{p-2}$ of ways of choosing a further $p-2$ points of $C$, times the number $e_1^{e_2-1}$ of ways of choosing a point from each of the $\C_1$ classes in $C^+\setminus C$ (to obtain a $(p+e_2)$-subset $Z'$ containing $X$ with the properties just discussed), times the number $r'$ of blocks containing $Z'$. Similarly, the number $\lambda_Y$ of blocks containing $Y$ is equal to the number $\binom{e_1-2}{p-1}$ of ways of choosing a further $p-1$ points of $C$, times the number $e_1^{e_2-2}$ of ways of choosing a point from each of the $\C_1$ classes in $C^+\setminus C$ apart from the $\C_1$-class containing $y_3$ (again to obtain a $(p+e_2)$-subset $Z'$ containing $Y$ with the properties just discussed), times the number $r'$ of blocks containing $Z'$. Thus 
    \[
    \lambda_X= \binom{e_1-3}{p-2}\cdot e_1^{e_2-1} \cdot r' \quad \hbox{and}\quad \lambda_Y= \binom{e_1-2}{p-1}\cdot e_1^{e_2-2} \cdot r'.
    \]
If $\D^s$ were a $3$-design then $\lambda_X=\lambda_Y$, and this holds if and only if $\binom{e_1-3}{p-2}\cdot e_1 = \binom{e_1-2}{p-1}$. A direct computation shows, in turn, that this is true if and only if $(p-1)e_1=e_2-2$, which is a contradiction. Thus $\D^s$ is not a $3$-design.

Now we show that $\D^s$ is a $2$-design.
By Corollary~\ref{cor:2des-red} and Remark~\ref{rem:2des-b}, we may prove $\D^s$ is a $2$-design by verifying condition \eqref{2des-b} in Remark \ref{rem:2des-b} for each $i \in \{1,2,\ldots, s-1\}$. 
Recall that for $i=1$ the condition \eqref{2des-b} is equivalent to \eqref{2des-1}, and note that
    \[ \frac{k(k-1)}{v-1} = \frac{(q+1)q}{q^2 + q} = 1. \]
Thus \eqref{2des-b} for $i \in \{1, 2, \ldots, s-1\}$ is equivalent to
    \begin{equation} \label{s>2}
        \begin{array}{lll}
        \sum_{(\delta_j)_{j>1} \in \prod_{j>1} \mathbb{Z}_{e_j}} x_{(\delta_j)_{j>1}} \left( x_{(\delta_j)_{j>1}} - 1 \right) &= e_1 - 1 &\text{for $i=1$, and} \\
        \sum_{(\delta_j)_{j>i-1} \in \prod_{j>i-1} \mathbb{Z}_{e_j}} x_{(\delta_j)_{j>i-1}} \left( x_{(\delta_j)_{j>i}} - x_{(\delta_j)_{j>i-1}} \right) &= (e_i - 1) \prod_{j=0}^{i-1} e_j  &\text{for $2 \leq i \leq s-1$}.
        \end{array}
    \end{equation}

For the first equation in \eqref{s>2} we use the values of $x_{(\delta_j)_{j>i-1}}$ and $x_{(\delta_j)_{j>i-1}} - 1$ in the second and fourth columns of Table \ref{tab:s>2}, replacing $i$ with $2$. Thus we have
    \begin{align*}
    \sum_{(\delta_j)_{j>1} \in \prod_{j>1} \mathbb{Z}_{e_j}} x_{(\delta_j)_{j>1}} \left( x_{(\delta_j)_{j>1}} - 1 \right)
    &= (p+1)p + \sum_{i=2}^s (e_i - 1) \cdot 1\cdot 0 \\
    &= p^2 + p \\
    &= e_1 - 1.
    \end{align*}

For $2 \leq i \leq s-1$, we use the values of $x_{(\delta_j)_{j>i-1}}$ and $x_{(\delta_j)_{j>i}} - x_{(\delta_j)_{j>i-1}}$ in the second  and fifth columns of Table \ref{tab:s>2}. This gives us (using Lemma \ref{lem:e_i}(i) in the last line)
    \begin{align*}
    &\sum_{(\delta_j)_{j>i-1} \in \prod_{j>i-1} \mathbb{Z}_{e_j}} x_{(\delta_j)_{j>i-1}} \left( x_{(\delta_j)_{j>i}} - x_{(\delta_j)_{j>i-1}} \right) \\
    &= \big( p^{2^{i-2}} + 1 \big)\big( p^{2^{i-1}} - p^{2^{i-2}} \big) + (e_i - 1) \cdot 1\cdot p^{2^{i-1}} + \sum_{i'=i+1}^s (e_{i'} - 1) \cdot 1\cdot 0 \\
    &= \big( p^{2^{i-2}} + 1 \big)(e_i - 1) + (e_i - 1)p^{2^{i-1}} \\
    &= (e_i - 1)\big( p^{2^{i-1}} + p^{2^{i-2}} + 1 \big) \\
    &= (e_i - 1) \prod_{j=1}^{i-1} e_j.
    \end{align*}

Therefore all equations in \eqref{s>2} hold, so $\D^s$ is a $2$-design by Corollary \ref{cor:2des-red}.

To compute $\lambda$, observe that $k(k-1) = (q+1)q = v-1$, so $\lambda = \frac{bk(k-1)}{v(v-1)} = \frac{b}{v}$ where $b = |\B|$. By Lemma \ref{lem:b}, $|\B| = \binom{e_1-1}{p} \cdot \frac{c_s}{p+1} \cdot \prod_{i=2}^s (c_{i-1})^{e_i-1}$; since $c_s = |\PP| = v$, the formula for $\lambda$ follows immediately.
\end{proof}

\section{Uniform case} \label{sec:uniform}

Before obtaining the constructions described in Section \ref{sec:ex}, the search for examples of block-transitive, $s$-chain-imprimitive $2$-designs was attempted by considering parameter sets that satisfy certain restrictions. One such restriction is the condition that the numbers $e_i$ are all equal to some constant $e$. In this case $v = e^s$, so by Corollary \ref{cor:2des-red}, the point-block structure $\D = \D(\underbrace{e, \ldots, e}_s; \,\chi)$ is a $2$-design exactly when
    \begin{equation} \label{2des-e-1}
    \sum_{C \in \C_1} x_C(x_C - 1) = \frac{k(k-1)}{e^s - 1} (e-1)
    \end{equation}
and
    \begin{equation} \label{2des-e}
    \sum_{C \in \C_{i-1}} x_{C} \left( x_{C^+} - x_{C} \right) = \frac{k(k-1)}{e^s-1} (e-1)e^{i-1} \quad \text{for $2 \leq i \leq s-1$},
    \end{equation}
where, for each $C \in \C_{i-1}$, $C^+$ is as defined in \eqref{C+}. Note that the last assertion of Theorem \ref{mainthm:2des} is equivalent to the condition that the quantities in \eqref{2des-e-1} and \eqref{2des-e} are even integers.

We found a few examples of $2$-designs $\D(e,e,e; \,\chi)$ and one example of a $2$-design $\D(e,e,e,e; \,\chi)$, which we describe in Examples~\ref{ex:e=3} and~\ref{ex:s=4}. As yet no general construction is known (see Question~\ref{Q4}).

\begin{example} \label{ex:e=3}
Assume that $s = 3$ and $e = 3$, so that $v = 3^3 = 27$. Take $\PP = \mathbb{Z}_3 \times \mathbb{Z}_3 \times \mathbb{Z}_3$, $G = S_3 \wr S_3 \wr S_3$, and $\B = B^G$ where
    \begin{align*}
    B &= \big\{ (0,0,0), (1,0,0), (0,1,0), (1,1,0), (0,2,0), (1,2,0), (0,0,1), (1,0,1), (0,1,1), (1,1,1), \\
    &\phantom{=}\;\; (0,0,2), (1,0,2), (0,1,2) \big\}.
    \end{align*}
Then $k = 13$. Recall that the $\C_1$-classes are the sets
    \[ 
    C_{(\delta_2,\delta_3)} = \big\{ (\varepsilon_1,\delta_2,\delta_3) \ \big| \  \varepsilon_1 \in \mathbb{Z}_3 \big\} 
    \]
for $(\delta_2,\delta_3) \in \mathbb{Z}_3 \times \mathbb{Z}_3$, and the $\C_2$-classes are the sets
    \[ 
    C_{(\delta_3)} = \big\{ (\varepsilon_1,\varepsilon_2,\delta_3) \ \big| \  (\varepsilon_1,\varepsilon_2) \in \mathbb{Z}_3 \times \mathbb{Z}_3 \big\} 
    \]
for $\delta_3 \in \mathbb{Z}_3$. The values of $x_{(\delta_2,\delta_3)} = \big| B \cap C_{(\delta_2,\delta_3)} \big|$ for $C_{(\delta_2,\delta_3)} \in \C_1(B)$, and of $x_{(\delta_3)} = \big| B \cap C_{(\delta_3)} \big|$ for $C_{(\delta_3)} \in \C_2(B)$, with $\C_1(B)$ and $\C_2(B)$ as defined in \eqref{C(X)}, are listed in Table \ref{tab:e=3}.

\begin{table}[ht]
    \centering
    \begin{tabular}{ccc}
    \hline
    $(\delta_2,\delta_3)$ & $x_{(\delta_2,\delta_3)}$ & $x_{(\delta_3)}$ \\
    \hline\hline
    $(0,0)$ & $2$ & $6$ \\
    $(1,0)$ & $2$ &  \\
    $(2,0)$ & $2$ &  \\
    \hline
    $(0,1)$ & $2$ & $4$ \\
    $(1,1)$ & $2$ &  \\
    \hline
    $(0,2)$ & $2$ & $3$ \\
    $(1,2)$ & $1$ &  \\
    \hline
    \end{tabular}
    \caption{Array of $B$ in Example \ref{ex:e=3}}
    \label{tab:e=3}
\end{table}

By Corollary~\ref{cor:2des-red} and Remark~\ref{rem:2des-b}, we may prove that $\D$ is a $2$-design by showing that 
    \begin{equation} \label{e=3}
    \begin{cases}
    &\displaystyle\sum_{(\delta_2,\delta_3) \in \mathbb{Z}_3 \times \mathbb{Z}_3} x_{(\delta_2,\delta_3)} \left( x_{(\delta_2,\delta_3)} - 1 \right) = \frac{k(k-1)}{e^3 - 1} (e-1), \ \text{and} \\
    &\displaystyle\sum_{(\delta_2,\delta_3) \in \mathbb{Z}_3 \times \mathbb{Z}_3} x_{(\delta_2,\delta_3)} \left(x_{(\delta_3)} - x_{(\delta_2,\delta_3)} \right) = \frac{k(k-1)}{e^3 - 1}(e-1)e.
    \end{cases}
    \end{equation}
Note that the first condition comes from \eqref{2des-e-1}, and the second is from \eqref{2des-e} with $i=2$. Observe that
    \[ \frac{k(k-1)}{e^3-1}(e-1) = \frac{13\cdot 12}{26} \cdot 2 = 12. \]
Hence, using the values in Table \ref{tab:e=3} we get
    \[ \sum_{(\delta_2,\delta_3) \in \mathbb{Z}_3 \times \mathbb{Z}_3} x_{(\delta_2,\delta_3)} \left( x_{(\delta_2,\delta_3)} - 1 \right) = 6 \cdot 2(2-1) = 12, \]
and
    \begin{align*}
    \sum_{(\delta_2,\delta_3) \in \mathbb{Z}_3 \times \mathbb{Z}_3} x_{(\delta_2,\delta_3)} \left( x_{(\delta_3)} - x_{(\delta_2,\delta_3)} \right)
    &= 3 \cdot 2(6-2) + 2 \cdot 2(4-2) + 2(3-2) + 1(3-1) \\
    &= 36 \\
    &= 12\cdot 3.
    \end{align*}
So the conditions \eqref{e=3} are both satisfied, which proves that $\D$ is a $2$-design.

The number $b = |\B|$ of point subsets whose array is equivalent to $\chi_B$ is $b = 2^2 \cdot 3^{10}$, so we have
    \[ \lambda = \frac{bk(k-1)}{v(v-1)} = \frac{3^{10} \cdot 2^2 \cdot 13 \cdot 12}{3^3 \cdot 13 \cdot 2} = 2^3 \cdot 3^8. \]
\end{example}

\begin{example} \label{ex:s=4}
Assume that $s = 4$ and $e = 3$, so that $v = 3^4 = 81$. Let $\PP = \mathbb{Z}_3 \times \mathbb{Z}_3 \times \mathbb{Z}_3 \times \mathbb{Z}_3$, $G = S_3 \wr S_3 \wr S_3 \wr S_3$, and $\B = B^G$ where
    \begin{align*}
    B &= \big\{ (0,0,0,0), (0,1,0,0), (0,2,0,0), (0,0,1,0), (1,0,1,0), (0,0,2,0), (1,0,2,0), (0,0,0,1), \\
    &\phantom{=}\;\; (0,1,0,1), (0,2,0,1), (0,0,1,1), (1,0,1,1), (0,0,2,1), (0,0,0,2), (0,1,0,2), (0,2,0,2) \big\}.
    \end{align*}
Then $k = 16$, and $G$ leaves invariant the partitions
    \begin{itemize}
        \item $\C_1$ with classes $C_{(\delta_2,\delta_3,\delta_4)} = \big\{ (\varepsilon_1,\delta_2,\delta_3,\delta_4) \ \big| \ \varepsilon_1 \in \mathbb{Z}_3 \big\}$ for each $(\delta_2,\delta_3,\delta_4) \in \mathbb{Z}_3 \times \mathbb{Z}_3 \times \mathbb{Z}_3$,
        \item $\C_2$ with classes $C_{(\delta_3,\delta_4)} = \big\{ (\varepsilon_1,\varepsilon_2,\delta_3,\delta_4) \ \big| \ (\varepsilon_1,\varepsilon_2) \in \mathbb{Z}_3 \times \mathbb{Z}_3 \big\}$ for each $(\delta_3,\delta_4) \in \mathbb{Z}_3 \times \mathbb{Z}_3$, and
        \item $\C_3$ with classes $C_{(\delta_4)} = \big\{ (\varepsilon_1,\varepsilon_2,\varepsilon_3,\delta_4) \ \big| \ (\varepsilon_1,\varepsilon_2,\varepsilon_3) \in \mathbb{Z}_3 \times \mathbb{Z}_3 \times \mathbb{Z}_3 \big\}$ for each $\delta_4 \in \mathbb{Z}_3$.
    \end{itemize}
The values of $x_{(\delta_2,\delta_3,\delta_4)}$, $x_{(\delta_3,\delta_4)}$, and $x_{(\delta_4)}$ corresponding to the classes in $\C_i(B)$ for $1 \leq i \leq 3$ are given in Table \ref{tab:s=4}.

\begin{table}[ht]
    \centering
    \begin{tabular}{cccc}
    \hline
    $(\delta_2,\delta_3,\delta_4)$ & $x_{(\delta_2,\delta_3,\delta_4)}$ & $x_{(\delta_3,\delta_4)}$ & $x_{(\delta_4)}$ \\
    \hline\hline
    $(0,0,0)$ & $1$ & $3$ & $7$ \\
    $(1,0,0)$ & $1$ &  &  \\
    $(2,0,0)$ & $1$ &  &  \\
    \hline
    $(0,1,0)$ & $2$ & $2$ & $7$ \\
    \hline
    $(0,2,0)$ & $2$ & $2$ & $7$ \\
    \hline\hline
    $(0,0,1)$ & $1$ & $3$ & $6$ \\
    $(1,0,1)$ & $1$ &  &  \\
    $(2,0,1)$ & $1$ &  &  \\
    \hline
    $(0,1,1)$ & $2$ & $2$ & $6$ \\
    \hline
    $(0,2,1)$ & $1$ & $1$ & $6$ \\
    \hline\hline
    $(0,0,2)$ & $1$ & $3$ & $3$ \\
    $(1,0,2)$ & $1$ &  &  \\
    $(2,0,2)$ & $1$ &  &  \\
    \hline
    \end{tabular}
    \caption{Array of $B$ in Example \ref{ex:s=4}}
    \label{tab:s=4}
\end{table}

By Corollary~\ref{cor:2des-red} and Remark~\ref{rem:2des-b}, we may prove that $\D$ is a $2$-design by showing that the following conditions  are satisfied. They are obtained from \eqref{2des-e} (and \eqref{2des-e-1} for $i=1$).
    \begin{equation} \label{e=4}
    \begin{cases}
    &\sum_{(\delta_2,\delta_3,\delta_4) \in \mathbb{Z}_3^3} x_{(\delta_2,\delta_3,\delta_4)} \big( x_{(\delta_2,\delta_3,\delta_4)} - 1 \big) = \frac{k(k-1)}{e^4-1} (e-1) \\
    &\sum_{(\delta_2,\delta_3,\delta_4) \in \mathbb{Z}_3^3} x_{(\delta_2,\delta_3,\delta_4)} \big( x_{(\delta_3,\delta_4)} - x_{(\delta_2,\delta_3,\delta_4)} \big) = \frac{k(k-1)}{e^4-1} e(e-1) \\
    &\sum_{(\delta_3,\delta_4) \in \mathbb{Z}_3^2} x_{(\delta_3,\delta_4)} \big( x_{(\delta_4)} - x_{(\delta_3,\delta_4)} \big) = \frac{k(k-1)}{e^4-1} e^2(e-1)
    \end{cases}
    \end{equation}
Note that
    \[ \frac{k(k-1)}{e^4-1} (e-1) = \frac{16\cdot 15}{81-1} \cdot 2 = 6. \]
To verify the first equation in \eqref{e=4} we use the values of $x_{(\delta_2,\delta_3,\delta_4)}$  listed in Table \ref{tab:s=4}. We get
    \begin{align*}
    \sum_{(\delta_2,\delta_3,\delta_4) \in \mathbb{Z}_3^3} x_{(\delta_2,\delta_3,\delta_4)} \big( x_{(\delta_2,\delta_3,\delta_4)} - 1 \big)
    &= 3 \cdot 2(2-1)
    = 6,
    \end{align*}
so the first equation in \eqref{e=4} holds. For the second equation in \eqref{e=4} we have
    \begin{align*}
    &\sum_{(\delta_2,\delta_3,\delta_4) \in \mathbb{Z}_3^3} x_{(\delta_2,\delta_3,\delta_4)} \big( x_{(\delta_3,\delta_4)} - x_{(\delta_2,\delta_3,\delta_4)} \big) \\
    &= 3 \cdot 1(3-1) + 2 \cdot 2(2-2) + 3 \cdot 1(3-1) + 2(2-2) + 1(1-1) + 3 \cdot 1(3-1) \\
    &= 6 + 0 + 6 + 0 + 0 + 6 \\
    &= 18 = 6\cdot 3,
    \end{align*}
which shows that the second equation is also true. Finally, for the third equation in \eqref{e=4} we have
    \begin{align*}
    \sum_{(\delta_3,\delta_4) \in \mathbb{Z}_3^2} x_{(\delta_3,\delta_4)} \big( x_{(\delta_4)} - x_{(\delta_3,\delta_4)} \big)
    &= 3(7-3) + 2 \cdot 2(7-2) + 3(6-3) + 2(6-2) + 1(6-1) +  3(3-3) \\
    &= 12 + 20 + 9 + 8 + 5 + 0 \\
    &= 54 = 6\cdot 3^2,
    \end{align*}
so likewise the third equation in \eqref{e=4} holds. Therefore $\D$ is indeed a $2$-design.

The number $b = |\B|$ of subsets of $\PP$ whose array is equivalent to $\chi_B$ is $b = 2^2 \cdot 3^{21}$, and therefore
    \[
    \lambda = \frac{bk(k-1)}{v(v-1)} = \frac{3^{21} \cdot 2^3 \cdot 2^4 \cdot 3 \cdot 5}{3^4 \cdot 2^4 \cdot 5} = 2^2 \cdot 3^{18}.
    \]
\end{example}

The remaining examples we found all have $s=3$, so $v=e^3$. The values of $e, k,$ and a generating block $B$ for these examples are listed in Table \ref{tab:e=s=3}. The proof that these produce $2$-designs, and the computation of $\lambda$, is similar to the arguments used in Example \ref{ex:e=3}.

\begin{table}[ht]
    \centering  \renewcommand{\arraystretch}{1.1}
    \begin{tabular}{rrrl}
    \hline
    $e$ & $k$ & $\lambda$ & $B$ \\
    \hline
    $4$ & $7$ & $2^{14} \cdot 3^2$ & $\big\{ (0,0,0), (1,0,0), (0,0,1), (0,1,1), (0,0,2), (0,1,2), (0,2,2) \big\}$ \\
    \hline
    $4$ & $15$ & $2^{24} \cdot 3^4 \cdot 5$ & \parbox[t]{11.2cm}{$\big\{ (0,0,0), (1,0,0), (2,0,0), (0,1,0), (0,2,0), (0,3,0), (0,0,1), (1,0,1),$ $(0,1,1), (1,1,1), (0,0,2), (0,1,2), (0,2,2), (0,0,3), (0,1,3) \big\}$} \\
    \hline
    $9$ & $14$ & $2^{11} \cdot 3^{37} \cdot 5 \cdot 7^2$ & \parbox[t]{11.2cm}{$\big\{ (0,0,0), (1,0,0), (0,0,1), (0,1,1), (0,0,2), (0,1,2), (0,0,3), (0,1,3),$
        $(0,0,4), (0,1,4), (0,2,4), (0,3,4), (0,0,5), (0,0,6) \big\}$} \\
    \hline
    \end{tabular}
    \caption{Blocks $B$ generating $2$-designs $\D(e,e,e; \,\chi_B)$}
    \label{tab:e=s=3}
\end{table}


\begin{thebibliography}{}
\bibitem{M2}
A. Montinaro, The symmetric  $2-(v,k.\lambda)$ designs, with $k > \lambda(\lambda-3)/2$, admitting a flag-transitive, point-imprimitive automorphism group are known.  ArXiv:  2212.08893.
\end{thebibliography}
\end{document}